\pdfoutput = 1
\documentclass[letterpaper]{article}
\usepackage{aaai_mod}
\usepackage{times}
\usepackage{helvet}
\usepackage{courier}
\usepackage{amsmath,amssymb,amsfonts}
\usepackage{algorithmic}
\usepackage{graphicx,epstopdf}
\usepackage{textcomp}
\usepackage{xcolor}
\usepackage{amsthm}
\usepackage{enumitem}
\usepackage{algorithm}
\usepackage{bbm} %indicator function
\usepackage{multicol}
\usepackage{multirow}
\usepackage{array}
\usepackage{cases}
\usepackage[subrefformat=parens, labelformat=parens]{subfig}
\allowdisplaybreaks

\usepackage{makecell} % for stacking texts

\newtheorem{theorem}{Theorem}
\newtheorem{definition}{Definition}
\newtheorem{lemma}{Lemma}

\newtheorem{assumption}{Assumption}

\makeatletter
\newcommand{\thickhline}{%
    \noalign {\ifnum 0=`}\fi \hrule height 1pt
    \futurelet \reserved@a \@xhline
}
\newcolumntype{"}{@{\hskip\tabcolsep\vrule width 1pt\hskip\tabcolsep}}
\makeatother

\begin{document}
\title{Riemannian stochastic recursive momentum method for non-convex optimization}
\author{Andi Han \and Junbin Gao\\
University of Sydney\\
\{andi.han, junbin.gao\}@sydney.edu.au
}
\maketitle
\begin{abstract}
We propose a stochastic recursive momentum method for Riemannian non-convex optimization that achieves a near-optimal complexity of $\Tilde{\mathcal{O}}(\epsilon^{-3})$ to find $\epsilon$-approximate solution with one sample. That is, our method requires $\mathcal{O}(1)$ gradient evaluations per iteration and does not require restarting with a large batch gradient, which is commonly used to obtain the faster rate. Extensive experiment results demonstrate the superiority of our proposed algorithm.
\end{abstract}

\section{Introduction}
We consider the problem of expectation (online) minimization over Riemannian manifold $\mathcal{M}$, defined as
\begin{equation}
    \min_{x \in \mathcal{M}} F(x) := \mathbb{E}_{\omega}[f(x, \omega)] \label{expectation_problem},
\end{equation}
where $F: \mathcal{M} \xrightarrow[]{} \mathbb{R}$, is a sufficiently smooth and potentially non-convex function. When $\omega$ can be finitely sampled from its support $\Omega$, problem \eqref{expectation_problem} reduces to empirical risk (finite-sum) minimization, $\min_{x \in \mathcal{M}} F(x) := \frac{1}{n} \sum_{i=1}^n f_i(x)$, where $n$ is the number of component functions. 

In this paper, we focus on the case where full gradient of ${F}(x)$ is inaccessible as in online setting or when $n$ is extremely large under finite-sum setting. We specifically consider optimization on Riemannian manifold, which includes Euclidean space as a special case. The Riemannian framework is ubiquitous in a variety of contexts. For example, in machine learning and statistics, principal component analysis (PCA) and matrix completion can be formulated on Grassmann manifold \cite{ZhangRSVRG2016,BoumalRTRMC2011}. In image processing, the tasks of diffusion tensor imaging segmentation and clustering can be cast on symmetric positive definite (SPD) manifold \cite{ChengDTIP2012}. Joint diagonalization for independent component analysis (ICA) is a problem over Stiefel manifold, which is useful for image separation \cite{FaridICA1999}.

Riemannian geometry provides the minimal ingredients that allow unconstrained optimization methods to be properly defined. A fundamental choice to solve problem \eqref{expectation_problem} is Riemannian stochastic gradient descent (RSGD) \cite{BonnabelRSGD2013}, which generalizes SGD on Euclidean space \cite{RobbinsSA1951}. Recently, SGD with coordinate-wise adaptive learning rates have become predominately popular within Deep Learning community. This includes Adagrad \cite{DuchiAdaGrad2011}, RMSProp \cite{HintonRMSPROP2012} and Adam \cite{kingmaADAM2014}, just to name a few. The main feature of these methods is to scale gradient adaptively so that each parameter has different learning rate. This has been proved to improve robustness and escape saddle points faster compared to standard SGD \cite{StaibAGESP2019}. For Riemannian optimization of problem \eqref{expectation_problem}, the lack of a canonical coordinate system and the highly nonlinear geometry make it difficult to extend such adaptation effectively. Regardless, \citeauthor{KumarCRMSPROP2018} (\citeyear{KumarCRMSPROP2018}) proposed constrained SGD with momentum (cSGD-M) and constrained RMSProp (cRMSProp) that adapt learning rates by coordinate-wise operations on matrix manifolds. However unlike Euclidean space, parallel transporting past gradients will likely distort gradient features, such as sparsity. Also, they fail to provide convergence guarantee. \citeauthor{BecigneulRADAM2018} (\citeyear{BecigneulRADAM2018}) generalized Adam-like adaptation and momentum to a product of Riemannian manifolds (referred to as RADAM and RAMSGRAD) where adaptation is implemented across component manifolds. \citeauthor{LiCayleyADAM2020} (\citeyear{LiCayleyADAM2020}) introduced Cayley Adam tailored for Stiefel manifold, exploiting its unique geometry. The only work that proves non-convex convergence on matrix manifolds is \cite{KasaiAdaptiveManifold2019}, where they propose RASA that adapts row and column subspaces of underlying manifold, ensuring a convergence rate of $\Tilde{\mathcal{O}}(1/\sqrt{T})$. This matches that of RSGD up to a logarithmic factor.

\begin{table*}[!t]
\setlength{\tabcolsep}{3pt}
\renewcommand{\arraystretch}{1.4}
\centering
\caption{Comparison of methods for Riemannian online non-convex optimization}
\begin{tabular}{c |ccccc}
\thickhline
             & Complexity & Large Batch & Small Batch & Restarting & Manifold types \\ \hline
\makecell{RSRG/RSPIDER \\ \cite{HanVRABA2020,ZhouRSPIDER2019}} &$\mathcal{O}(\epsilon^{-3})$            & $\mathcal{O}(\epsilon^{-2})$             & $\mathcal{O}(\epsilon^{-1})$             & Yes  & General       \\
\makecell{RSGD \\ \cite{HosseiniRSGD2017}}         & $\mathcal{O}(\epsilon^{-4})$            & ---             & $\mathcal{O}(1)$             & No  & General        \\
\makecell{RASA \\ \cite{KasaiAdaptiveManifold2019}}         & $\Tilde{\mathcal{O}}(\epsilon^{-4})$           & ---            & $\mathcal{O}(1)$            & No    & Matrix manifolds     \\
RSRM (this work)         &  $\Tilde{\mathcal{O}}(\epsilon^{-3})$          & $\mathcal{O}(1)$             & $\mathcal{O}(1)$             & No  & General  \\
\thickhline
\end{tabular}
\end{table*}

Although SGD-based methods enjoy low sampling cost, i.e. $\mathcal{O}(1)$ per iteration (one-sample), the main bottleneck of this method, similar to its Euclidean counterpart, is the unvanished gradient variance. This largely slows down its convergence. For this problem, variance reduction (VR) techniques are gaining increasing attentions. Many methods, including Riemannian stochastic variance reduced gradient (RSVRG) \cite{ZhangRSVRG2016,SatoRSVRG2019}, stochastic recursive gradient (RSRG) \cite{KasaiRSRG2018}, stochastic path-integrated differential estimator (RSPIDER) \cite{ZhouRSPIDER2019,ZhangRSPIDER2018} are generalized from their Euclidean versions. The main idea is to correct for stochastic gradient deviation by periodically computing a large batch gradient. As a result, gradient variance decreases as training progresses, which leads to a linear convergence. However, VR methods are originally designed for finite-sum optimization where full gradient access is possible. Some recent studies have extended Riemannian VR to online setting \cite{ZhouRSPIDER2019,HanVRABA2020}. Specifically, RSVRG requires $\mathcal{O}(\epsilon^{-10/3})$ stochastic gradient oracles to achieve $\epsilon$-approximate solution (see Definition \ref{definition_1}), which improves on $\mathcal{O}(\epsilon^{-4})$ of RSGD. Also, RSRG and RSPIDER require an even lower complexity of $\mathcal{O}(\epsilon^{-3})$. This rate has been proved to be optimal under one additional mean-squared smoothness assumption on Euclidean space \cite{ArjevaniLBSG2019}. 

Nevertheless, these online VR methods still require computing a large batch gradient, i.e. $\mathcal{O}(\epsilon^{-2})$ for each epoch and the mini-batch size for each inner iteration should also be at least $\mathcal{O}(\epsilon^{-1})$. To address this issue, we introduce a novel online variance reduction method, inspired by a recently proposed recursive momentum estimator \cite{CutkoskyMomentumVR2019}. Our main \textbf{contributions} are summarized below. 
\begin{itemize}
    \item We propose a Riemannian stochastic recursive momentum (RSRM) method that achieves a gradient complexity of $\Tilde{\mathcal{O}}(\epsilon^{-3})$ for online non-convex optimization, matching the lower bound up to a logarithmic factor.
    \item RSRM requires only $\mathcal{O}(1)$ gradient computations per iteration and does not need restarting with a large batch gradient. Thus, our method preserves the efficiency of SGD while achieving fast convergence as VR methods.
    \item Our convergence result holds for general manifolds while other online adaptive methods apply to restricted manifold types, such as matrix manifolds \cite{KumarCRMSPROP2018,KasaiAdaptiveManifold2019}, product manifolds \cite{BecigneulRADAM2018} and Stiefel manifold \cite{LiCayleyADAM2020}.
    \item Our extensive experiments prove that our method significantly outperforms other one-sample methods.
\end{itemize}
The rest of this paper is organized as follows. Section \ref{preliminary_section} introduces some useful definitions and notations as well as assumptions used for convergence analysis. Section \ref{algorithm_section} describes our proposed algorithms and highlights its relationships with RSGD, variance reduction and stochastic momentum. Section \ref{convergence_section} presents convergence analysis for RSRM and Section \ref{experiment_section} evaluates the proposed method on a variety of tasks and manifolds.

\section{Preliminaries}
\label{preliminary_section}
Riemannian manifold is a manifold with a smooth inner product $\langle \cdot, \cdot \rangle_x: T_x\mathcal{M} \times T_x\mathcal{M} \xrightarrow[]{} \mathbb{R}$ defined on tangent space $T_x\mathcal{M}$ for every $x \in \mathcal{M}$. It induces a norm on $T_x\mathcal{M}$: $\| u \|_x := \sqrt{\langle u, u \rangle_x}$. Retraction $R_x: T_x\mathcal{M} \xrightarrow{} \mathcal{M}$ maps a tangent vector to manifold surface satisfying some local conditions. That is, $R_x(0) = x$ and $\text{D}R_x(0)[u] = u$. The retraction curve is defined as $c(t) := R_x(t\xi)$ for $\xi \in T_x\mathcal{M}$. Denote $y = R_x(\xi)$. Then vector transport $\mathcal{T}_x^y$ (or equivalently $\mathcal{T}_\xi$) with respect to retraction $R$ maps $u \in T_x\mathcal{M}$ to $\mathcal{T}_x^y u \in T_x\mathcal{M}$ along the defined retraction curve $c(t)$. Note that exponential map $\text{Exp}_x$ is a special instance of retraction by restricting retraction curve to be a geodesic. Similarly, as a special case of vector transport, parallel transport $P_x^y$ maps a tangent vector in `parallel' along a curve while preserving its norm and direction. In this paper, we consider the more general and computationally efficient retraction and vector transport. Therefore our results can be trivially applied to exponential map and parallel transport. Implicitly, we consider only isometric vector transport $\mathcal{T}_x^y$, which satisfies $\langle u , v \rangle_x = \langle \mathcal{T}_x^y u , \mathcal{T}_x^y v \rangle_y$ for all $u, v \in T_x\mathcal{M}$. 

\textbf{Notations.} For the discussion that follows, we omit the subscript for norm and inner product, which should be clear from the context. We define a sampling set $\mathcal{S} = \{ \omega_{1},...,\omega_{|\mathcal{S}|} \}$ with cardinality $|\mathcal{S}|$. Each $\omega_{(\cdot)}$ is sampled independently from $\Omega$. We thus denote the Riemannian stochastic gradient $\text{grad} f_{\mathcal{S}}(x) := \frac{1}{|\mathcal{S}|} \sum_{\omega \in \mathcal{S}} \text{grad}f(x, \omega) \in T_x\mathcal{M}$. We denote $g(t) = \mathcal{O}(h(t))$ if there exists a positive constant $M$ and $t_0$ such that $g(t) \leq M h(t)$ for all $t \geq t_0$. We use $\Tilde{\mathcal{O}}(\cdot)$ to hide poly-logarithmic factors. We generally refer to $\| \cdot \|$ as the induced norm on tangent space of Riemannian manifold and use $\| \cdot \|_F$ to represent matrix Frobenius norm. Now we are ready to make some assumptions as follows. 

\begin{assumption}
\label{assump_1}
Iterates generated by RSRM stay continuously in a neighbourhood $\mathcal{X} \subseteq \mathcal{M}$ that contains an optimal point $x^*$.
The objective $F$ is continuously differentiable and has bounded initial suboptimality. That is, for all $x \in \mathcal{X}$, $F(x_1) - F(x^*) \leq \Delta$.
\end{assumption}

\begin{assumption}
\label{BiasVarAssump}
 Stochastic gradient $\text{grad}f(x, \omega)$ is unbiased with bounded variance. That is, for all $x \in \mathcal{X}, \omega \in \Omega$, it satisfies that
 \begin{align*}
     &\mathbb{E}_\omega \emph{grad}f(x, \omega) = \emph{grad}F(x), \\ &\mathbb{E}_\omega\| \emph{grad}f(x, \omega) - \emph{grad}F(x) \|^2 \leq \sigma^2,
 \end{align*}
 for some $\sigma > 0$.
\end{assumption}

\begin{assumption}
\label{retraction_L_assump}
 The objective $F$ is retraction $L$ smooth with respect to retraction $R$. That is, there exists a positive constant $L$ such that for all $x, y = R_x(\xi) \in \mathcal{X}$, we have
 \begin{equation*}
     F(y) \leq F(x) + \langle \emph{grad}F(x) ,\xi \rangle + \frac{L}{2} \| \xi \|^2.
 \end{equation*}
\end{assumption}

These three assumptions are standard in Riemannian stochastic gradient methods \cite{HosseiniRSGD2017,KasaiRSRG2018}. Note that the assumption of bounded iterates in neighbourhood $\mathcal{X}$ can be made with respect to the entire manifold $\mathcal{M}$, which results in stricter conditions on the retraction and vector transport in the following assumptions. To ensure retraction $L$ smoothness as in Assumption \ref{retraction_L_assump}, we require an upper-bounded Hessian property on the pullback function $f \circ R: T_x\mathcal{M} \xrightarrow[]{} \mathbb{R}$. That is, for all $x \in \mathcal{X}$ and $u \in T_x\mathcal{M}$ with unit norm, $\frac{d^2 f(R_x(tu))}{dt^2} \leq L$. In the work of RASA \cite{KasaiAdaptiveManifold2019}, the variance bound in Assumption \ref{BiasVarAssump} is replaced by $G$-gradient Lipschitz, which is $\| \text{grad}f(x, \omega) \| \leq G$. This however, amounts to a stronger assumption.

According to \cite{ArjevaniLBSG2019}, under the first three assumptions, SGD is minimax optimal. To obtain faster convergence, one further assumption of mean-squared retraction Lipschitz is required. This assumption is a straightforward generalization of mean-squared Lipschitz on Euclidean space, which is the minimal additional requirement to achieve the complexity lower bound. 

\begin{assumption}
\label{MSLipAssump}
 The objective $f$ is mean-squared retraction $\Tilde{L}$ Lipschitz. That is, there exists a positive constant $\Tilde{L}$ such that for all $x,y = R_x(\xi) \in \mathcal{X}, \omega \in \Omega$, 
 \begin{equation*}
     \mathbb{E}_\omega\| \emph{grad}f(x,\omega) - \mathcal{T}_y^x \emph{grad}f(y, \omega)  \|^2 \leq \Tilde{L}^2 \| \xi \|^2
 \end{equation*}
 holds with vector transport $\mathcal{T}_x^y$ along the retraction curve $c(t):= R_x(t\xi)$. 
\end{assumption}

It is aware that the standard assumption of retraction Lipschitzness is made with respect to parallel transport and one additional assumption that bounds the difference between vector transport and parallel transport is needed for Assumption \ref{MSLipAssump} to hold. See Lemma 4 in \cite{HanVRABA2020}. In this work, algorithm complexity is measured by the total number of stochastic first order oracles to achieve $\epsilon$-approximate solution, defined as follows. 

\begin{definition}[$\epsilon$-approximate solution and SFO]
\label{definition_1}
$\epsilon$-approximate solution by a stochastic algorithm is the output $x$ with expected square norm of its gradient less than $\epsilon^2$, i.e. $\mathbb{E}\| \emph{grad}F(x) \|^2 \leq \epsilon^2$. One stochastic first-order oracle (SFO) outputs a stochastic gradient $\emph{grad}f(x, \omega)$ given inputs $x$ and $\omega$ drawn from $\Omega$ \cite{GhadimiSFO2013}. 
\end{definition}

\section{Algorithms}
\label{algorithm_section}
\subsection{RSGD, variance reduction and momentum}
Riemannian stochastic gradient makes the following retraction update: $x_{t+1} = R_{x_t}( - \eta_t \text{grad}f_{\mathcal{S}_t}(x_t))$. This allows updates to follow the negative gradient direction while staying on the manifold. Variance reduction techniques utilize past gradient information to construct a modified gradient estimator with decreasing variance. In particular, the recursive gradient estimator in RSRG/RSPIDER achieves the optimal rate of $\mathcal{O}(\epsilon^{-3})$. That is, for each outer loop, a large batch gradient is computed as $d_0 = \text{grad}f_{\mathcal{S}_0}(x_0)$, where $|\mathcal{S}_0|$ is set to $n$ under finite-sum setting and $\mathcal{O}(\epsilon^{-2})$ under online setting. Within each inner iteration, stochastic gradient is corrected recursively based on its previous iterate:
\begin{equation}
    d_t = \text{grad}f_{\mathcal{S}_t} (x_t) - \mathcal{T}_{x_{t-1}}^{x_t} ( \text{grad}f_{\mathcal{S}_t}(x_{t-1}) - d_{t-1}), \label{recursive_estimator}
\end{equation}
where vector transport $\mathcal{T}_{x_{t-1}}^{x_t}$ is necessary to relate gradients on disjoint tangent spaces. To achieve the optimal complexity, mini-batch size $|\mathcal{S}_t|$ is set to be at least $\mathcal{O}(\epsilon^{-1})$. This choice of batch size can become very large, especially when we desire more accurate solutions.

On the other hand, stochastic gradient with momentum is not new on Euclidean space \cite{QianSGDM1999}, while the first paper that presents such an idea on Riemannian manifold is \cite{KumarCRMSPROP2018}. They simply takes a combination of current stochastic gradient and transported momentum, given as 
\begin{equation}
    d_t = \rho_t \mathcal{T}_{x_{t-1}}^{x_t} d_{t-1} + (1-\rho_t) \text{grad}f_{\mathcal{S}_t} (x_{t}), \label{momentum_estimator}
\end{equation}
where $\rho_t$ is commonly set to be $0.9$. However, no convergence analysis is provided. This idea has then been used in generalizing Adam and AMSGrad to Riemannian optimization \cite{BecigneulRADAM2018}, where they only established convergence on a product of manifolds for geodesically convex functions. Even on Euclidean space, the effectiveness of stochastic momentum over vanilla SGD has remained an open question.

\begin{algorithm}[!t]
 \caption{Riemannian SRM}
 \label{RSRM_algorithm}
 \begin{algorithmic}[1]
  \STATE \textbf{Input:} Step size $\eta_t$, recursive momentum parameter $\rho_t$, Initial point $x_1$.
  \STATE Compute $d_1 = \text{grad}f_{\mathcal{S}_1}(x_1)$.
  \FOR{$t = 1,...,T$}
  \STATE Update $x_{t+1} = R_{x_t} ( - \eta_t d_t)$.
  \STATE Compute $d_{t+1} = \text{grad}f_{\mathcal{S}_{t+1}}(x_{t+1}) + (1 - \rho_{t+1}) \mathcal{T}_{x_{t}}^{x_{t+1}} ( d_{t} - \text{grad}f_{\mathcal{S}_{t+1}}(x_{t}) )$.
  \ENDFOR
  \STATE \textbf{Output:} $\Tilde{x}$ uniformly chosen at random from $\{x_t \}_{t=1}^{T}$. 
 \end{algorithmic} 
\end{algorithm}

\subsection{Proposed RSRM}
Our proposed RSRM is presented in Algorithm \ref{RSRM_algorithm} where we make use of the recursive momentum estimator, originally introduced in \cite{CutkoskyMomentumVR2019,TranHSGD2019}:
\begin{align}
    d_{t} &= \rho_t \text{grad}f_{\mathcal{S}_t} (x_t) + (1 - \rho_t) \big( \text{grad}f_{\mathcal{S}_t} (x_t) \nonumber\\
    &- \mathcal{T}_{x_{t-1}}^{x_t} ( \text{grad}f_{\mathcal{S}_t}(x_{t-1}) - d_{t-1}) \big) \label{intermediate_d}\\
    &= \text{grad}f_{\mathcal{S}_t} (x_t) + (1- \rho_t)\mathcal{T}_{x_{t-1}}^{x_t} \big( d_{t-1} - \text{grad}f_{\mathcal{S}_t}(x_{t-1}) \big),
\end{align}
which hybrids stochastic gradient $\text{grad}f_{\mathcal{S}_t}(x_t)$ with the recursive gradient estimator in \eqref{recursive_estimator} for $\rho_t \in [0,1]$. This can be also viewed as combining momentum estimator in \eqref{momentum_estimator} with a scaled difference of $\text{grad}f_{\mathcal{S}_t}(x_t) - \mathcal{T}_{x_{t-1}}^{x_t} \text{grad}f_{\mathcal{S}_t}(x_{t-1})$. Note that we recover vanilla RSGD when $\rho_t = 0$ and the recursive estimator in \eqref{recursive_estimator} when $\rho_t = 1$. As we will demonstrate in Section \ref{convergence_section}, $\rho_t$ should be decreasing rather than fixed, thereby enabling a smooth transition from RSGD to RSRG. As a result, we do not require restarting the algorithm to achieve the optimal convergence. 

Compared with algorithm designs in Euclidean versions of SRM \cite{CutkoskyMomentumVR2019,TranHSGD2019}, our formulation and parameter settings are largely different. Specifically, \citeauthor{CutkoskyMomentumVR2019} (\citeyear{CutkoskyMomentumVR2019}) further adapts the recursive momentum parameter $\rho_t$ to the learning rate $\eta_t$ where the latter itself is adapted to the norm of stochastic gradient. This is claimed to relieve the parameter tuning process. However, they reintroduce three parameters, which are even less intuitive to be tuned (even though some are fixed to a default value). As shown in Section \ref{convergence_section}, we only require tuning the initial step size $\eta_0$ and initial momentum parameter $\rho_0$ (where the latter can be fixed to a good default value). Furthermore, the adaptive step size requires a uniform gradient Lipschitz condition, the same as in \cite{KasaiAdaptiveManifold2019} and also a uniform smoothness assumption, which is stronger than mean-sqaured smoothness in our setting. On the other hand, \citeauthor{TranHSGD2019} (\citeyear{TranHSGD2019}) replaces $\text{grad}f_{\mathcal{S}_t}(x_t)$ in \eqref{intermediate_d} with $\text{grad}f_{\mathcal{B}_t}(x_t)$ where $\mathcal{B}_t$ is independent with $\mathcal{S}_t$. This increases sampling complexity per iteration and also complicates its convergence analysis. In addition, they still require a large initial batch size $|\mathcal{S}_0| = \mathcal{O}(\epsilon^{-1})$ while our $|\mathcal{S}_0| = \mathcal{O}(1)$.

\section{Convergence results}
\label{convergence_section}
In this section, we prove convergence of RSRM. Define an increasing sigma-algebra $\mathcal{F}_t := \{ \mathcal{S}_1, ..., \mathcal{S}_{t-1} \}$. Hence by update rules in Algorithm \ref{RSRM_algorithm}, $x_t$ and $d_{t-1}$ are measurable in $\mathcal{F}_t$. We first present a Lemma that bounds the estimation error of the recursive momentum estimator. 

\begin{lemma}[Estimation error bound]
Suppose Assumptions \ref{assump_1} to \ref{MSLipAssump} hold and consider Algorithm \ref{RSRM_algorithm}. Then we can bound the expected estimation error of the estimator as
\label{variance_bound_estimator_lemma}
\begin{align}
    &\mathbb{E}\| d_{t+1} - \emph{grad}F(x_{t+1}) \|^2 \nonumber\\
    &\leq (1-\rho_{t+1})^2 \big( 1 + \frac{4\eta_{t}^2 \Tilde{L}^2}{|\mathcal{S}_{t+1}|} \big) \mathbb{E}\| d_{t} - \emph{grad}F(x_{t}) \|^2 \nonumber\\
    &+ \frac{4(1-\rho_{t+1})^2 \eta_{t}^2 \Tilde{L}^2}{|\mathcal{S}_{t+1}|} \mathbb{E}\| \emph{grad}F(x_{t}) \|^2 + \frac{2\rho_{t+1}^2 \sigma^2}{|\mathcal{S}_{t+1}|}.
\end{align}

\end{lemma}

The proof of this Lemma can be found in Supplementary material where it follows an idea similar to the bound in RSRG/RSPIDER \cite{HanVRABA2020}. The key difference is that we further use $\mathbb{E}\| d_t \|^2 \leq 2 \mathbb{E}\| d_t - \text{grad}F(x_t) \|^2 + 2 \mathbb{E}\| \text{grad}F(x_t) \|^2$ to show dependence on the full gradient. Based on the claims in \cite{CutkoskyMomentumVR2019}, we consider $\rho_t = \mathcal{O}(t^{-2/3})$ and $\eta_t = \mathcal{O}(t^{-1/3})$, so that $\mathbb{E}\| d_{t+1} - \text{grad}F(x_{t+1})\|^2 = \mathcal{O}(t^{-2/3} + \| \text{grad}F(x_t) \|^2)$. To see this, denote $s_{t+1} = d_{t+1} - \text{grad}F(x_{t+1})$. Then by noting that $(1 - \rho_t)^2 \leq 1 - \rho_t \leq 1$ and $\eta_t^2 \leq \eta_t \leq 1$,  Lemma \ref{variance_bound_estimator_lemma} suggests that $\mathbb{E}\| s_{t+1} \|^2 \leq \mathcal{O}( 1 - t^{-2/3}) \mathbb{E}\| s_t \|^2 + \mathcal{O}(t^{-2/3}) \mathbb{E}\| \text{grad}F(x_t) \|^2 + \mathcal{O}(t^{-4/3})$. Simply setting $\mathbb{E}\| s_{t+1} \|^2 = \mathbb{E}\| s_t \|^2$ yields the result. This implies that $\mathbb{E}\| \text{grad}F(x_t) \|^2 = \mathcal{O}(T^{-2/3})$, which matches the optimal rate of convergence. This claim is stated formally in the following Theorem. For simplicity, we consider $|\mathcal{S}_t| = b$ for all $t$.

\begin{theorem}[Convergence and complexity of RSRM]
\label{theorem_1}
Suppose Assumptions \ref{assump_1} to \ref{MSLipAssump} hold and consider Algorithm \ref{RSRM_algorithm} with $\eta_t = c_\eta (t+1)^{-1/3}$, $\rho_t = c_\rho t^{-2/3}$ where $c_\eta \leq \frac{1}{L}$ and $c_\rho = (\frac{10\Tilde{L}^2}{b} + \frac{1}{3}) c_\eta^2$. Then we have 
\begin{align*}
    \mathbb{E}\| \emph{grad}F(\Tilde{x}) \|^2 &= \frac{1}{T} \sum_{t=1}^T \mathbb{E} \| \emph{grad}F(x_t) \|^2 \\
    &\leq \mathcal{O}(\frac{M}{T^{2/3}}) = \Tilde{\mathcal{O}} (\frac{1}{T^{2/3}}),
\end{align*}
where $M := (6 \Delta + \frac{\sigma^2}{2\Tilde{L}^2} + \frac{\sigma^2 \ln(T+1)}{\Tilde{L}^2})/c_\eta$. To achieve $\epsilon$-approximate solution, we require an SFO complexity of $\Tilde{\mathcal{O}}(\epsilon^{-3})$. 
\end{theorem}
Proof of Theorem \ref{theorem_1} is included in Supplementary material. The proof idea is similar to \cite{CutkoskyMomentumVR2019} where we construct a Lyapunov function $R_t := \mathbb{E}[F(x_t)] + \frac{C}{\eta_{t-1}} \mathbb{E} \| d_t - \text{grad}F(x_t) \|^2 $ with $C = \frac{b}{12 \Tilde{L}^2}$. Theorem \ref{theorem_1} claims that RSRM achieves a near-optimal complexity of $\Tilde{\mathcal{O}}(\epsilon^{-3})$ with one-sample batch, i.e. $b = \mathcal{O}(1)$. And specifically under noiseless case where $\sigma^2 = 0$, we can further improve this result to the lower bound complexity $\mathcal{O}(\epsilon^{-3})$. One final remark can be made that our step size decays at a rate of $\mathcal{O}(t^{-1/3})$, which is slower compared to the SGD-based rate of $\mathcal{O}(t^{-1/2})$. This step size sequence is crucial for achieving the faster convergence, coupled with gradually reduced variance.

\section{Experiments}

\begin{figure*}[!ht]
\captionsetup{justification=centering}
    \centering
    \subfloat[Optimality gap vs. SFO (\texttt{SYN1}) \label{sfo_syn1}]{\includegraphics[width = 0.32\textwidth, height = 0.25\textwidth]{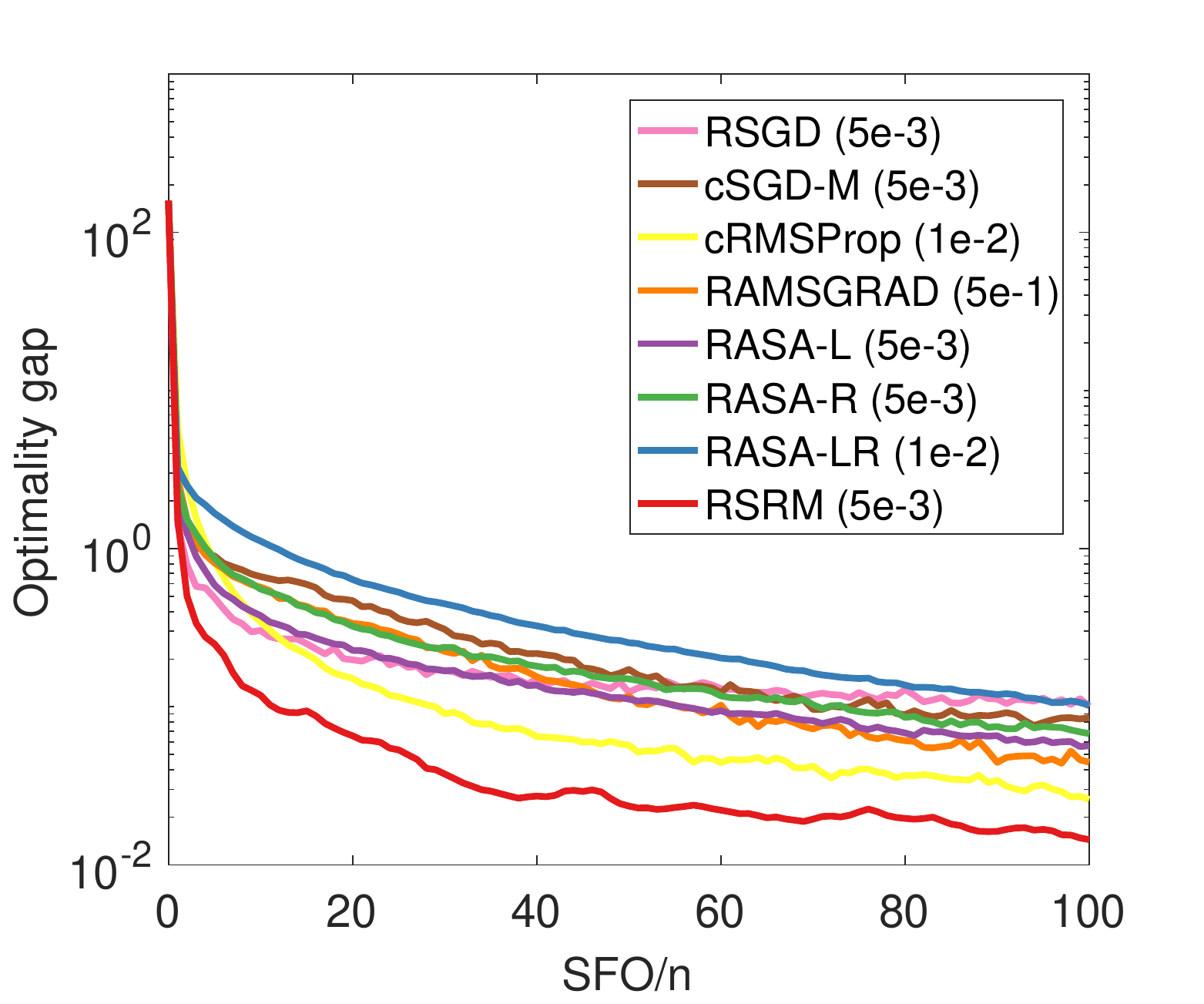}} 
    \hspace{0.01in}
    \subfloat[Optimality gap vs. SFO (\texttt{SYN2}) \label{syo_syn2}]{\includegraphics[width = 0.32\textwidth, height = 0.25\textwidth]{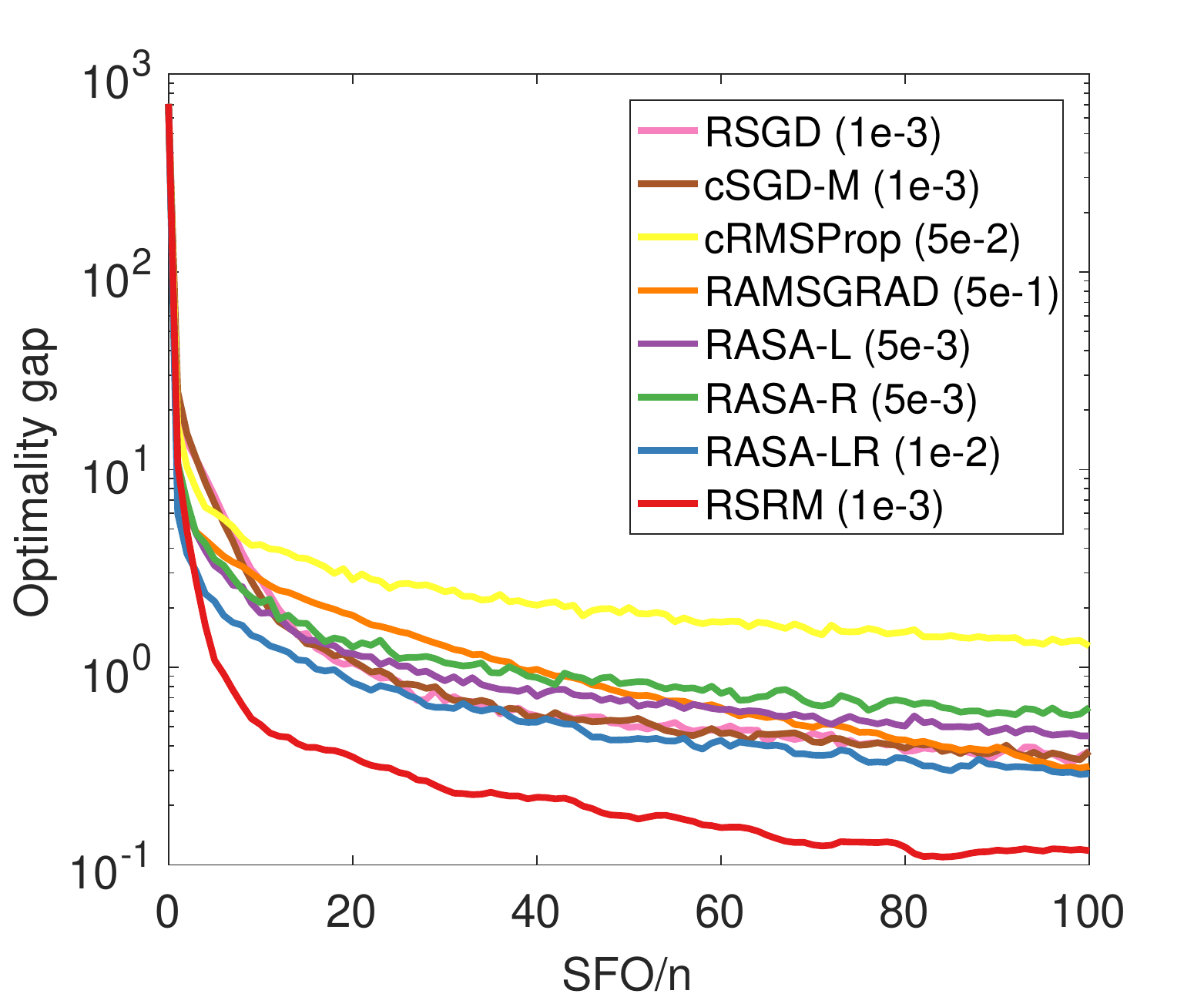}}
    \hspace{0.01in}
    \subfloat[Optimality gap vs. SFO (\texttt{SYN3}) \label{sfo_syn3}]{\includegraphics[width = 0.32\textwidth, height = 0.25\textwidth]{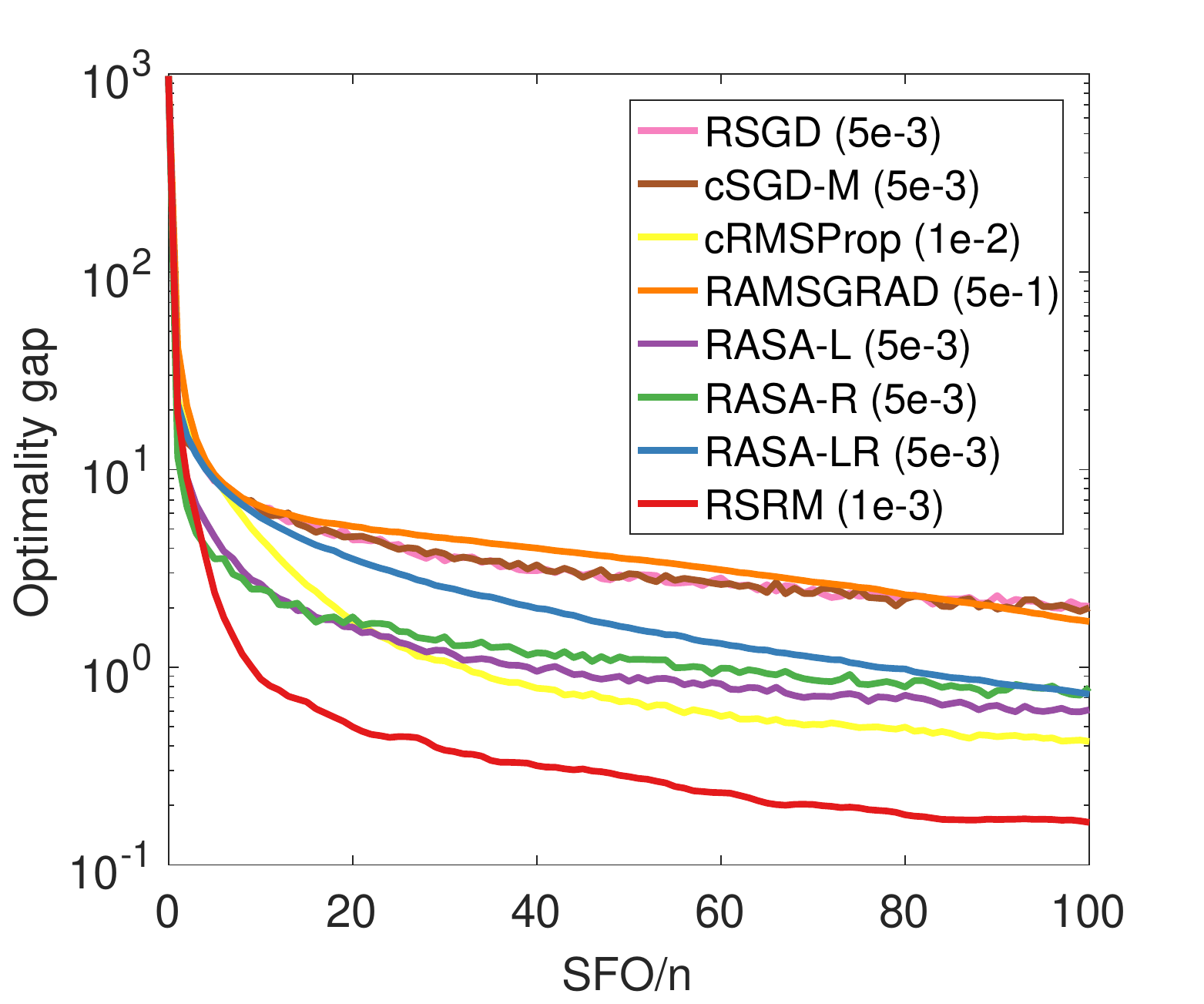}}\\[-0.2ex]
    \subfloat[Optimality gap vs. Runtime (\texttt{SYN1}) \label{time_syn1}]{\includegraphics[width = 0.32\textwidth, height = 0.25\textwidth]{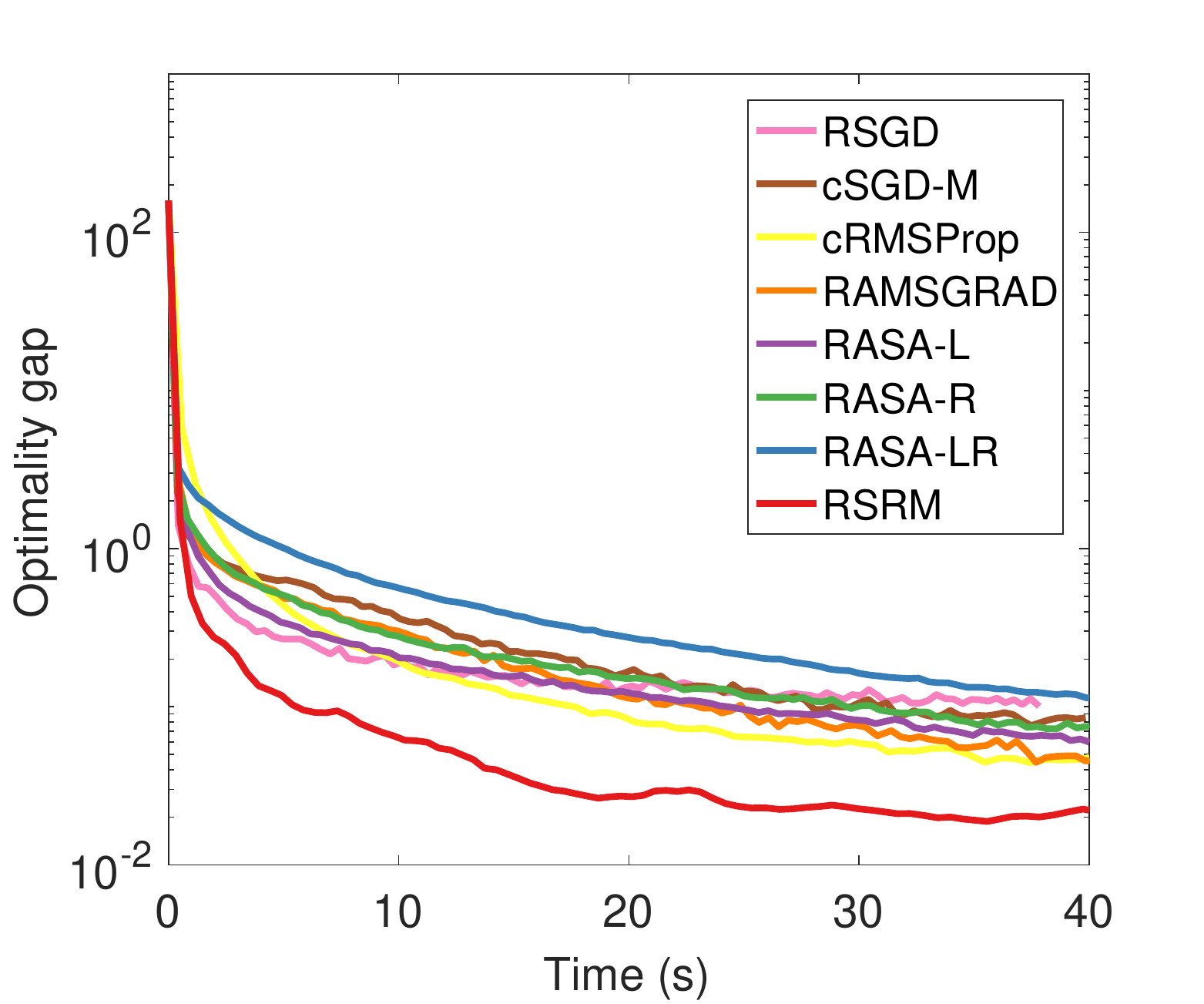}} 
    \hspace{0.01in}
    \subfloat[Optimality gap vs. SFO (\texttt{MNIST}) \label{sfo_mnist}]{\includegraphics[width = 0.32\textwidth, height = 0.25\textwidth]{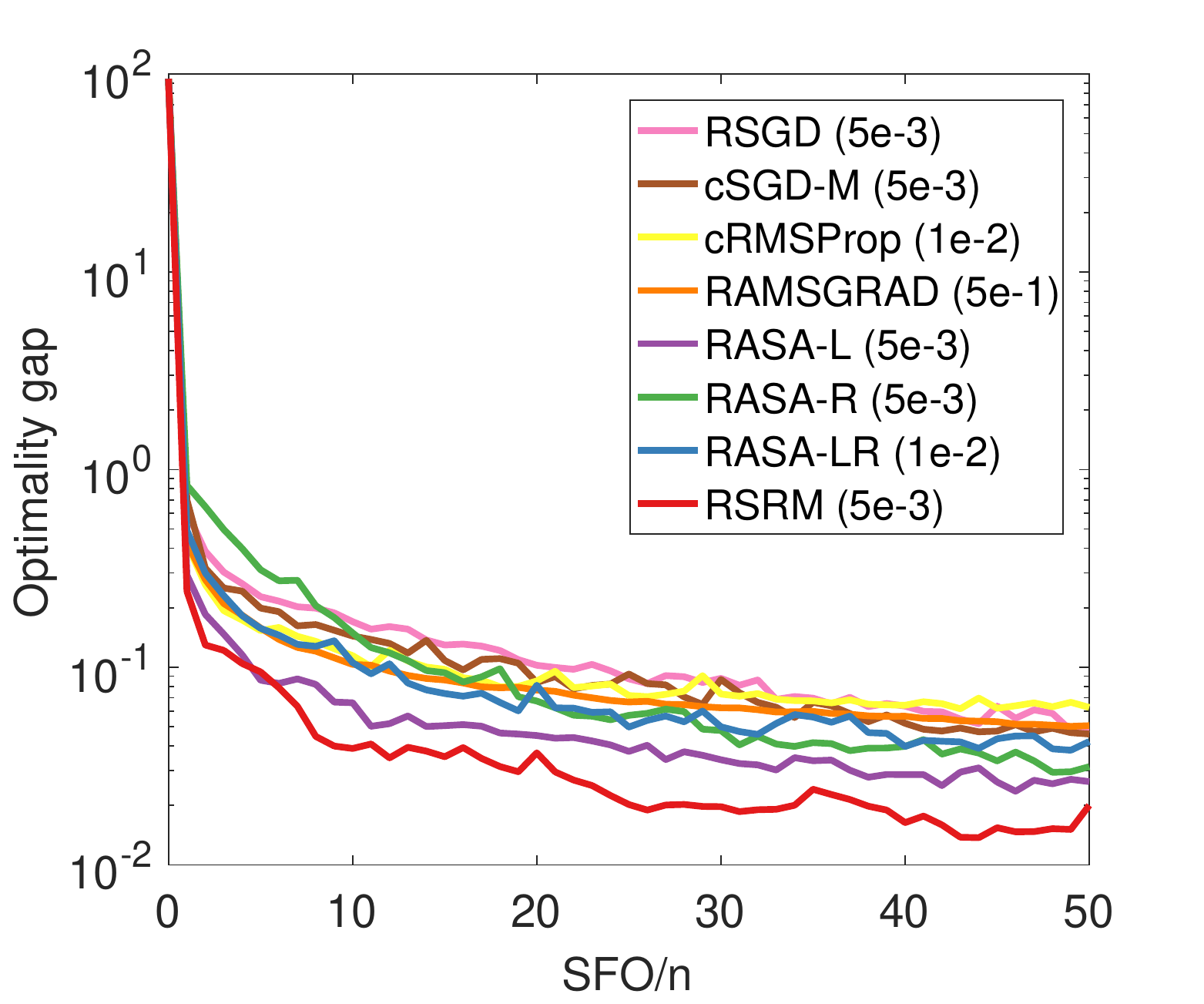}} 
    \hspace{0.01in}
    \subfloat[Optimality gap vs. SFO (\texttt{COVTYPE}) \label{sfo_covtype}]{\includegraphics[width = 0.32\textwidth, height = 0.25\textwidth]{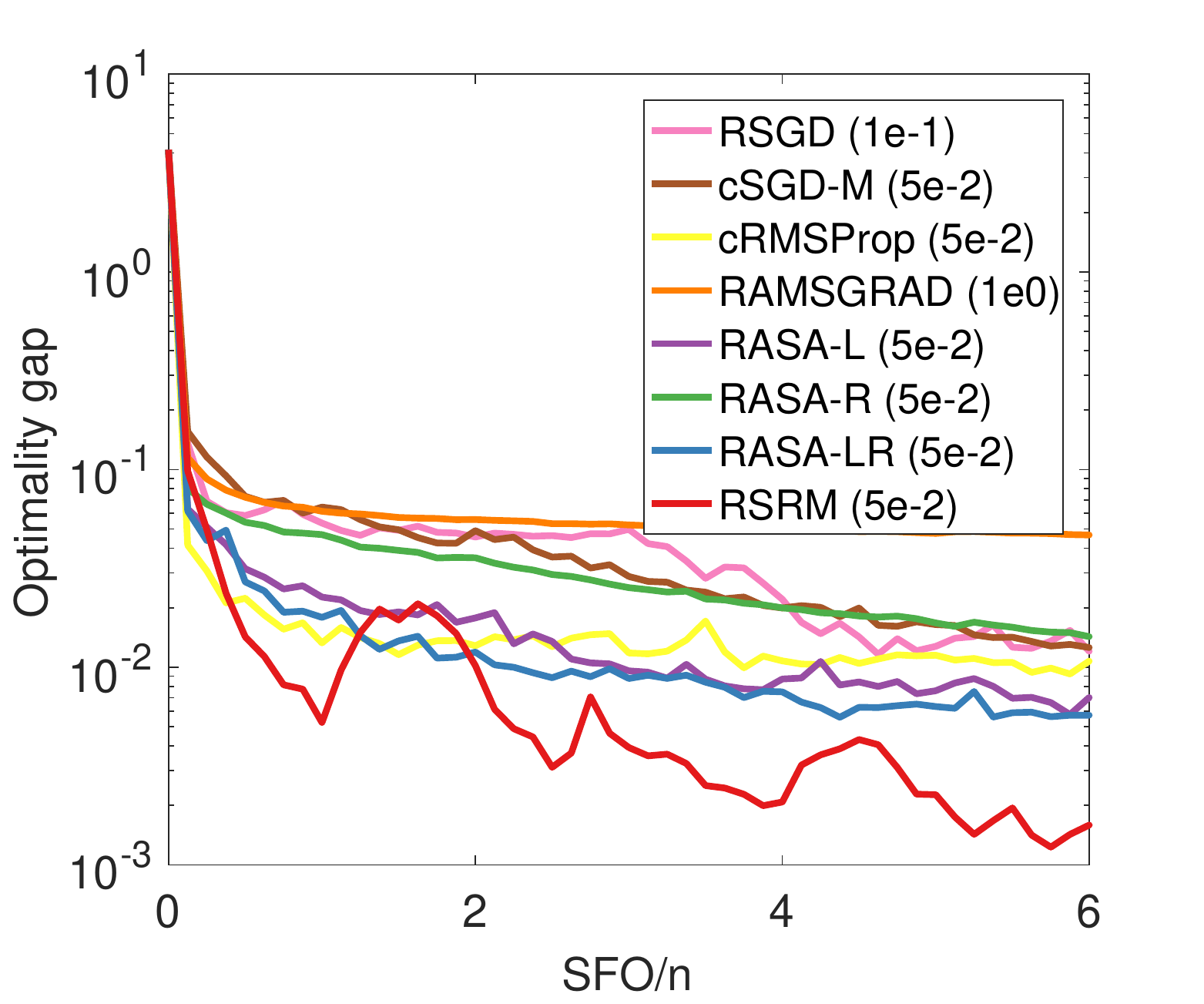}} \\[-0.2ex]
    \subfloat[Performance under different $\eta_0$ (\texttt{SYN1}) \label{eta_sensitivity}]{\includegraphics[width = 0.32\textwidth, height = 0.25\textwidth]{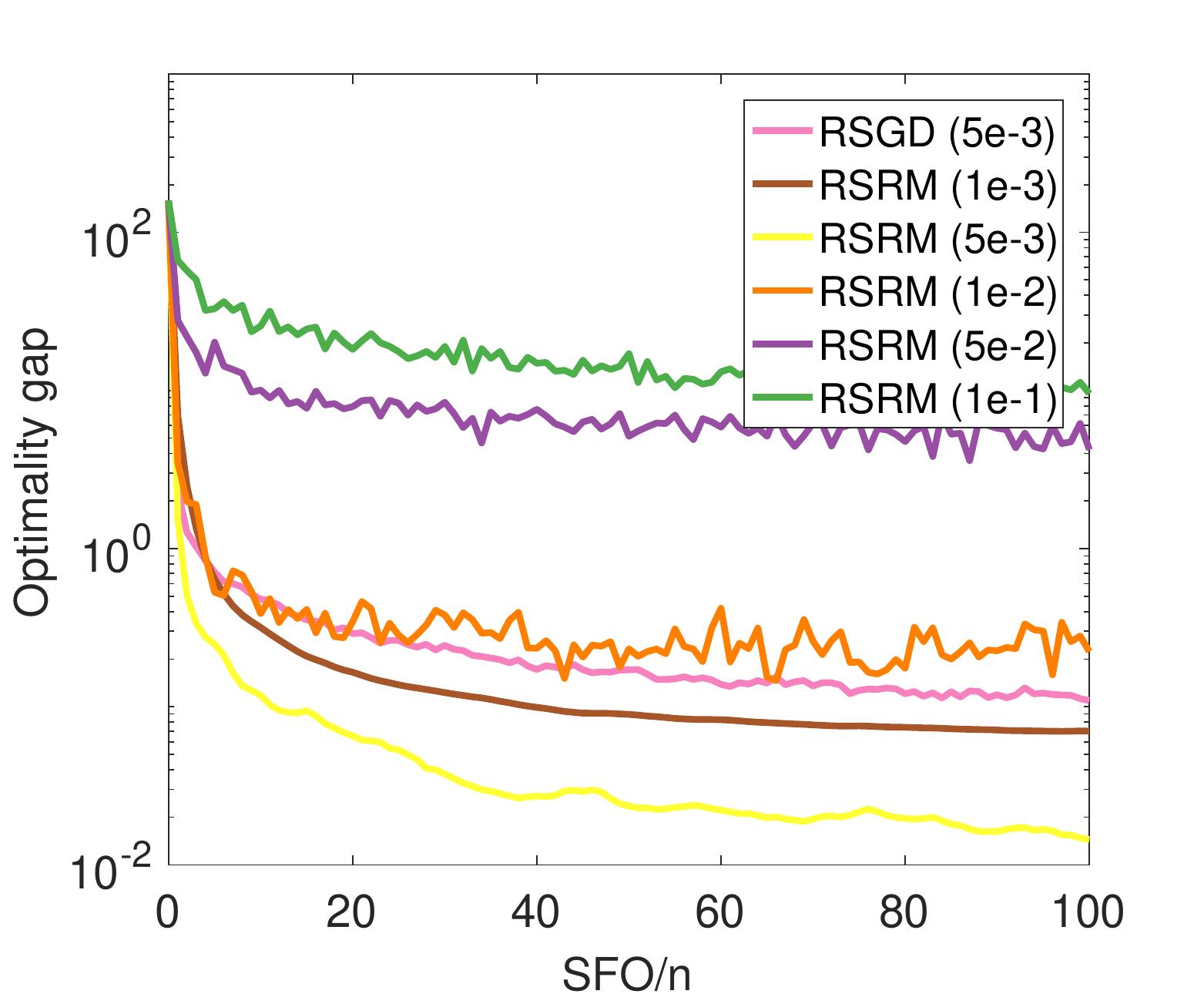}} 
    \hspace{0.01in}
    \subfloat[Performance under different $|\mathcal{S}_0|$ (\texttt{SYN1}) \label{b_sensitivity}]{\includegraphics[width = 0.32\textwidth, height = 0.25\textwidth]{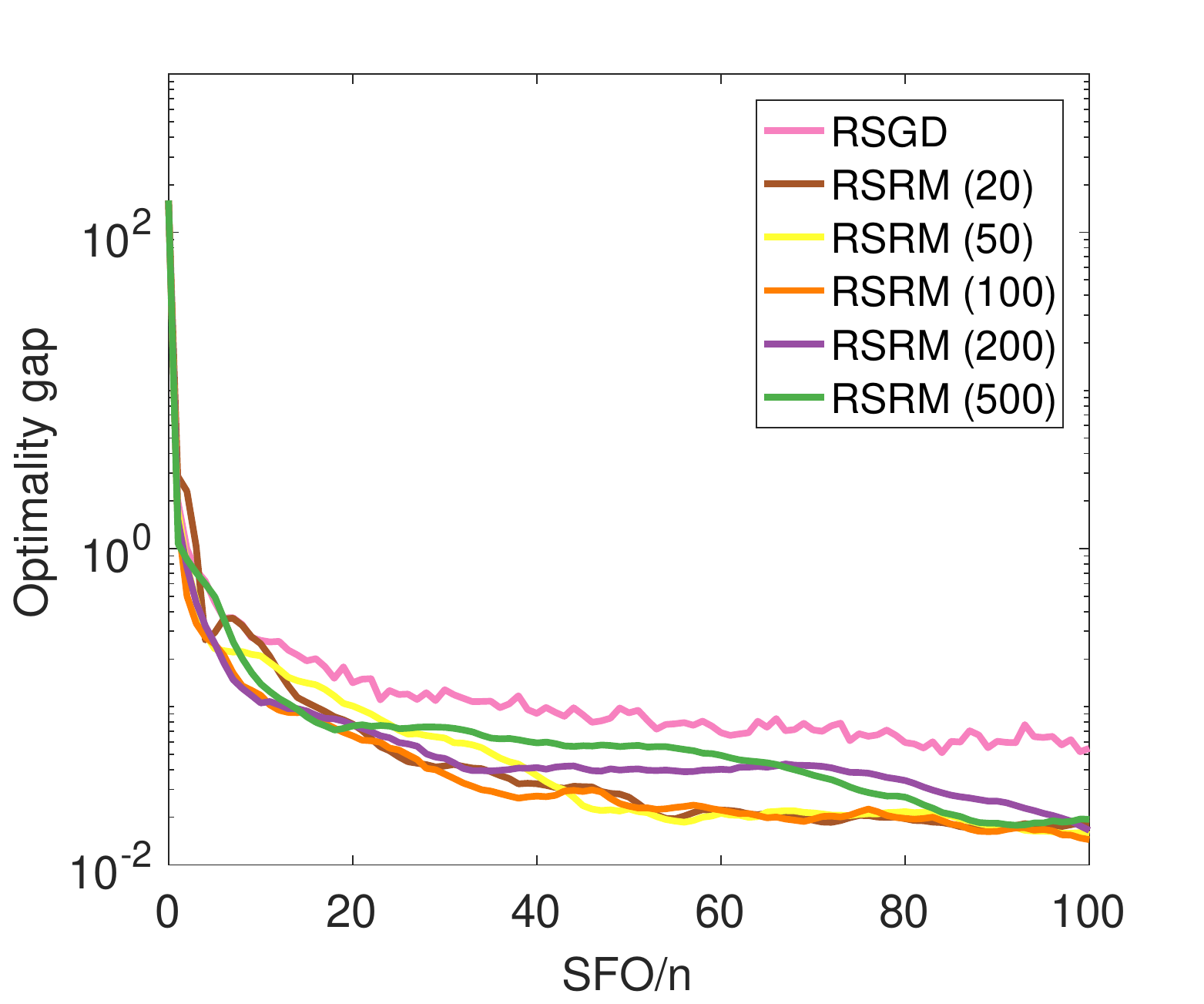}} 
    \caption{PCA problems on Grassmann manifold}
    \label{pca_grassmann}
\end{figure*}

\label{experiment_section}
In this section, we compare our proposed RSRM with other one-sample online methods, which are described as follows. The benchmark is the standard stochastic gradient method (RSGD) \cite{BonnabelRSGD2013}. We also consider cSGD-M and cRMSProp \cite{KumarCRMSPROP2018} where past gradients are transported by vector transport operator. For cRMSProp, we do not project and vector-transport its adaptation term, which is an element-wise square of stochastic gradient. Instead, we treat it as an element in the ambient Euclidean space and therefore we only project the resulting scaled gradient after applying this term. This modification turns out to yield significantly better convergence compared to its original design. Also we compare with RAMSGRAD \cite{BecigneulRADAM2018}, which is proposed for a product of manifolds. We thus modify the gradient momentum similar as in \cite{KumarCRMSPROP2018} while accumulating square norm of gradient instead of element-wise square. Hence, it only adapts the step size rather than the gradient. Finally, we consider RASA \cite{KasaiAdaptiveManifold2019} that adapts column and row subspaces of matrix manifolds. We similarly label its variants as RASA-L, RASA-R and RASA-LR to respectively represent adapting row (left) subspace, column (right) subspace and both. 

\begin{figure*}[!ht]
\captionsetup{justification=centering}
    \centering
    \subfloat[Optimality gap vs. SFO (\texttt{YALEB}) \label{sfo_yaleb_ica}]{\includegraphics[width = 0.32\textwidth, height = 0.25\textwidth]{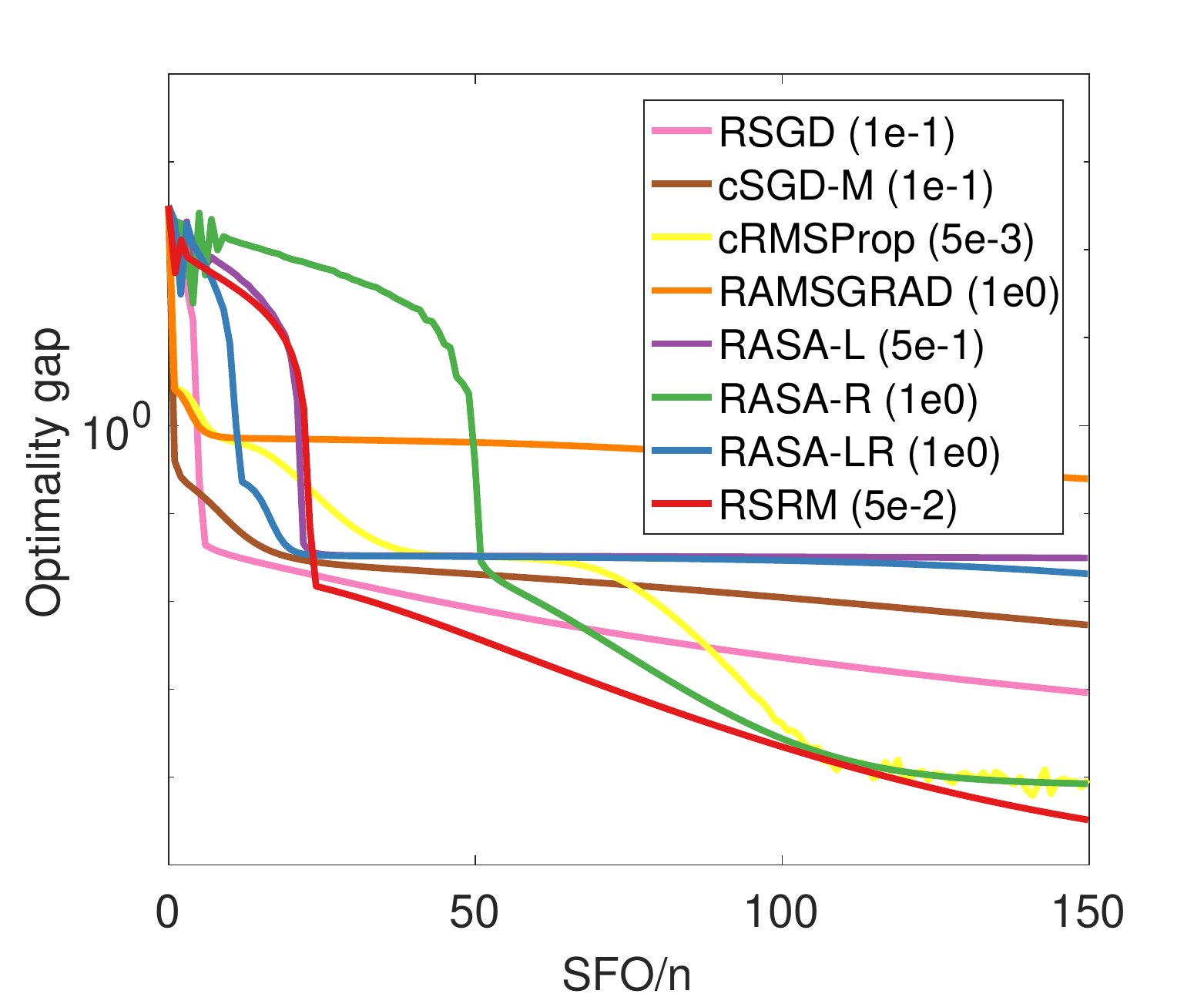}} 
    \hspace{0.01in}
    \subfloat[Optimality gap vs. SFO (\texttt{CIFAR100})]{\includegraphics[width = 0.32\textwidth, height = 0.25\textwidth]{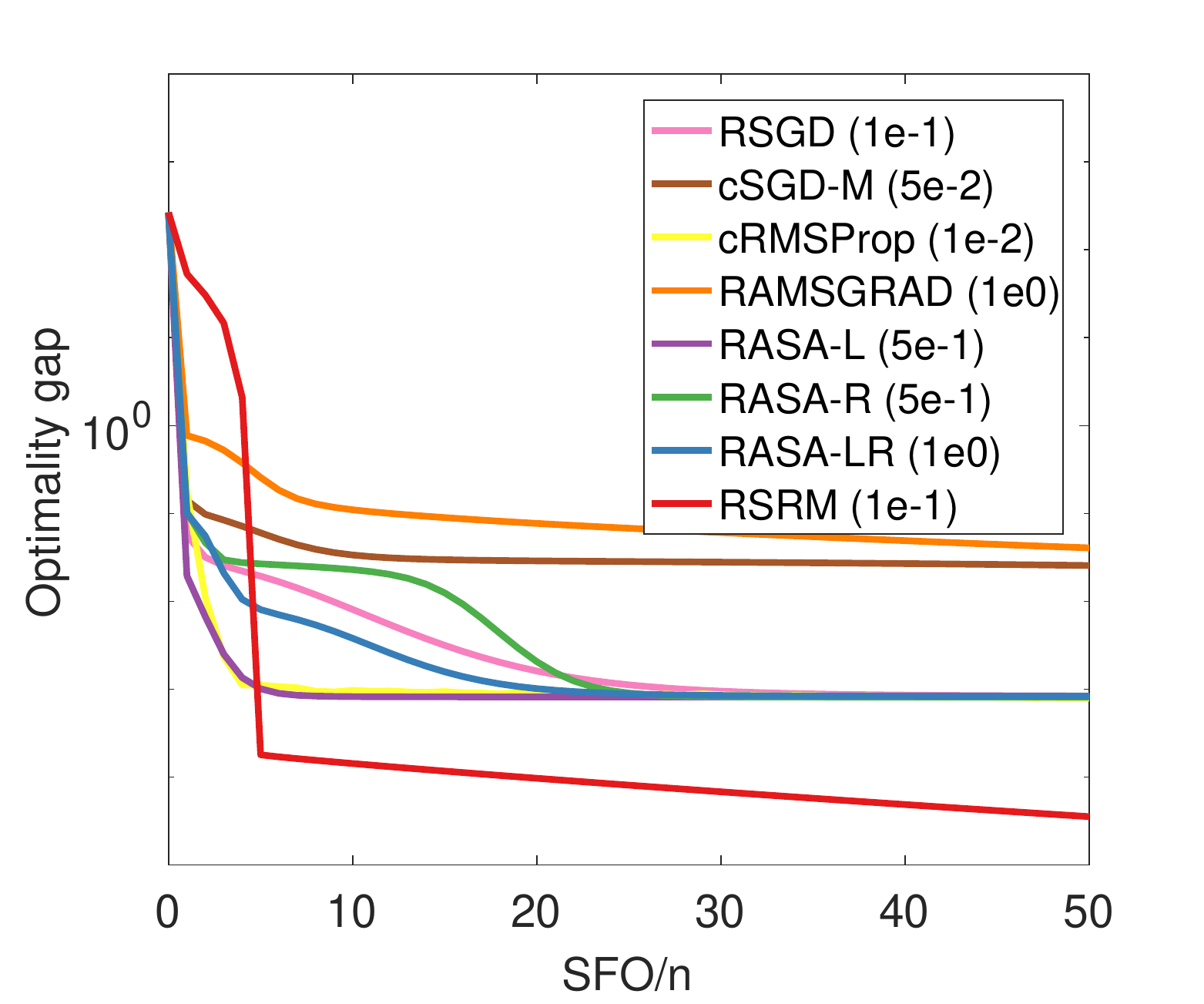}}
    \hspace{0.01in}
    \subfloat[Optimality gap vs. SFO (\texttt{COIL100})]{\includegraphics[width = 0.32\textwidth, height = 0.25\textwidth]{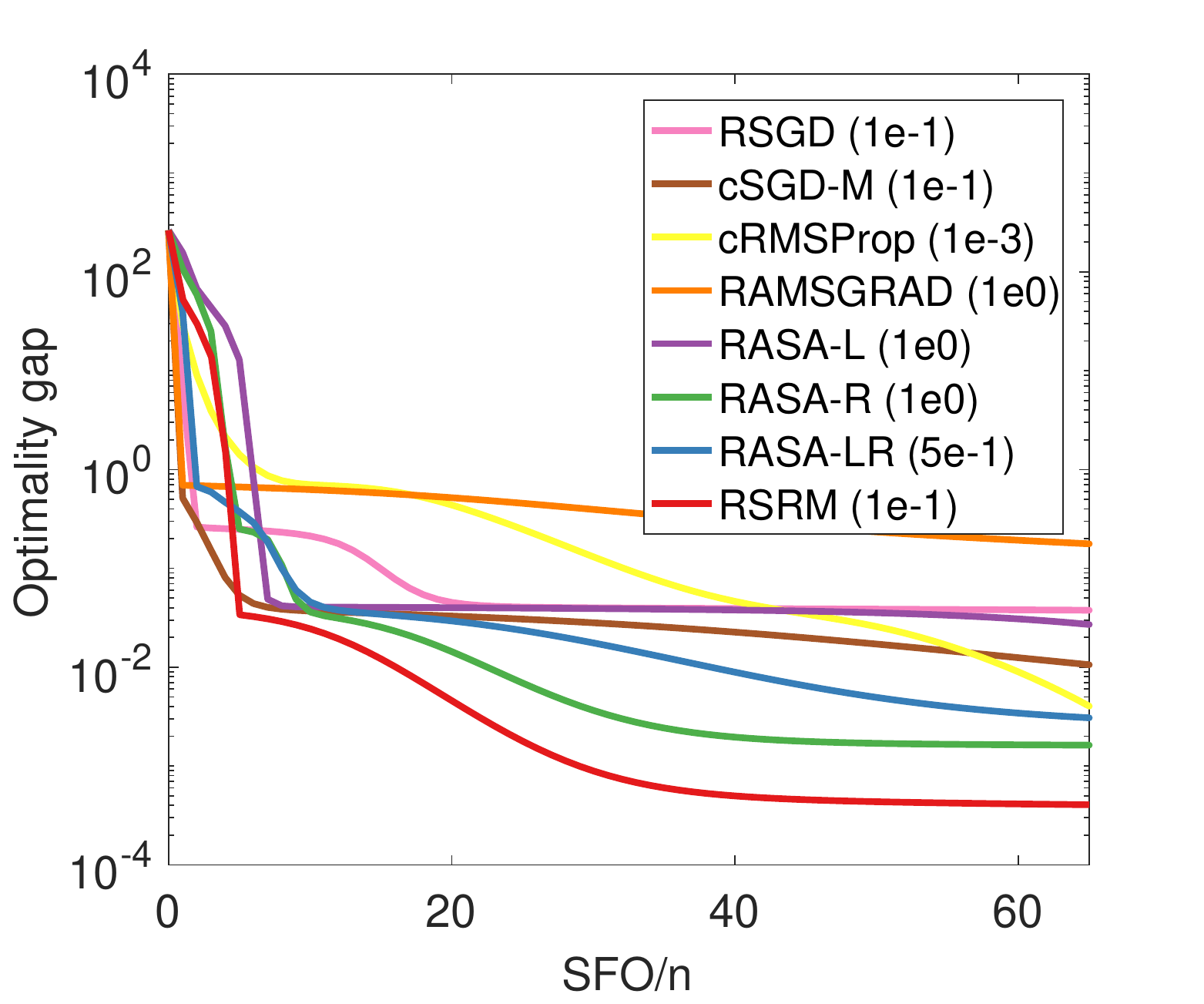}}\\
    \caption{ICA problems on Stiefel manifold}
    \label{ica_stiefel}
\end{figure*}

\begin{figure*}[!ht]
\captionsetup{justification=centering}
    \centering
    \subfloat[Optimality gap vs. SFO (\texttt{SYN1}) \label{sfo_syn_rc}]{\includegraphics[width = 0.32\textwidth, height = 0.25\textwidth]{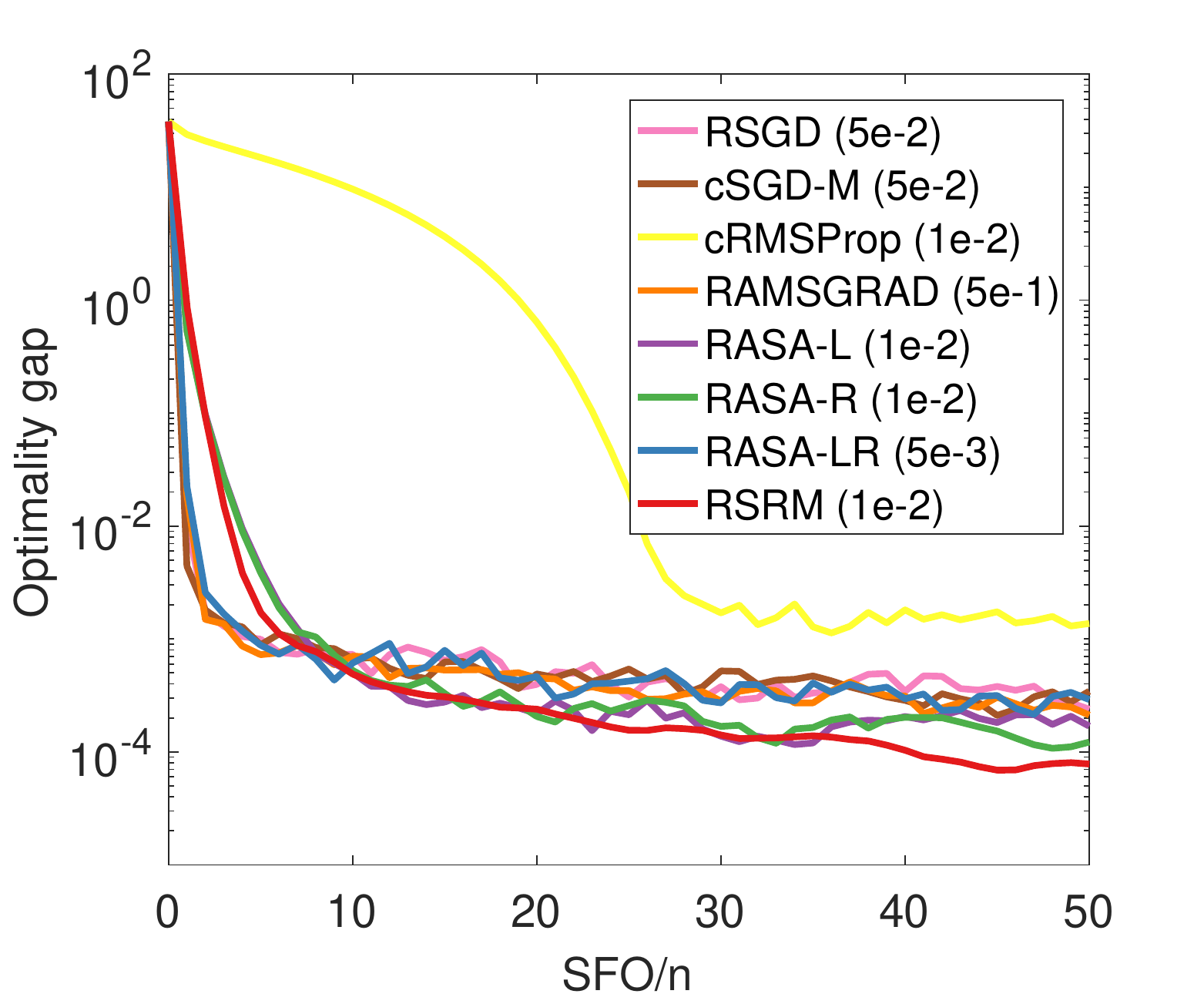}} 
    \hspace{0.01in}
    \subfloat[Optimality gap vs. SFO (\texttt{YABLEB}) \label{sfo_yaleb_rc}]{\includegraphics[width = 0.32\textwidth, height = 0.25\textwidth]{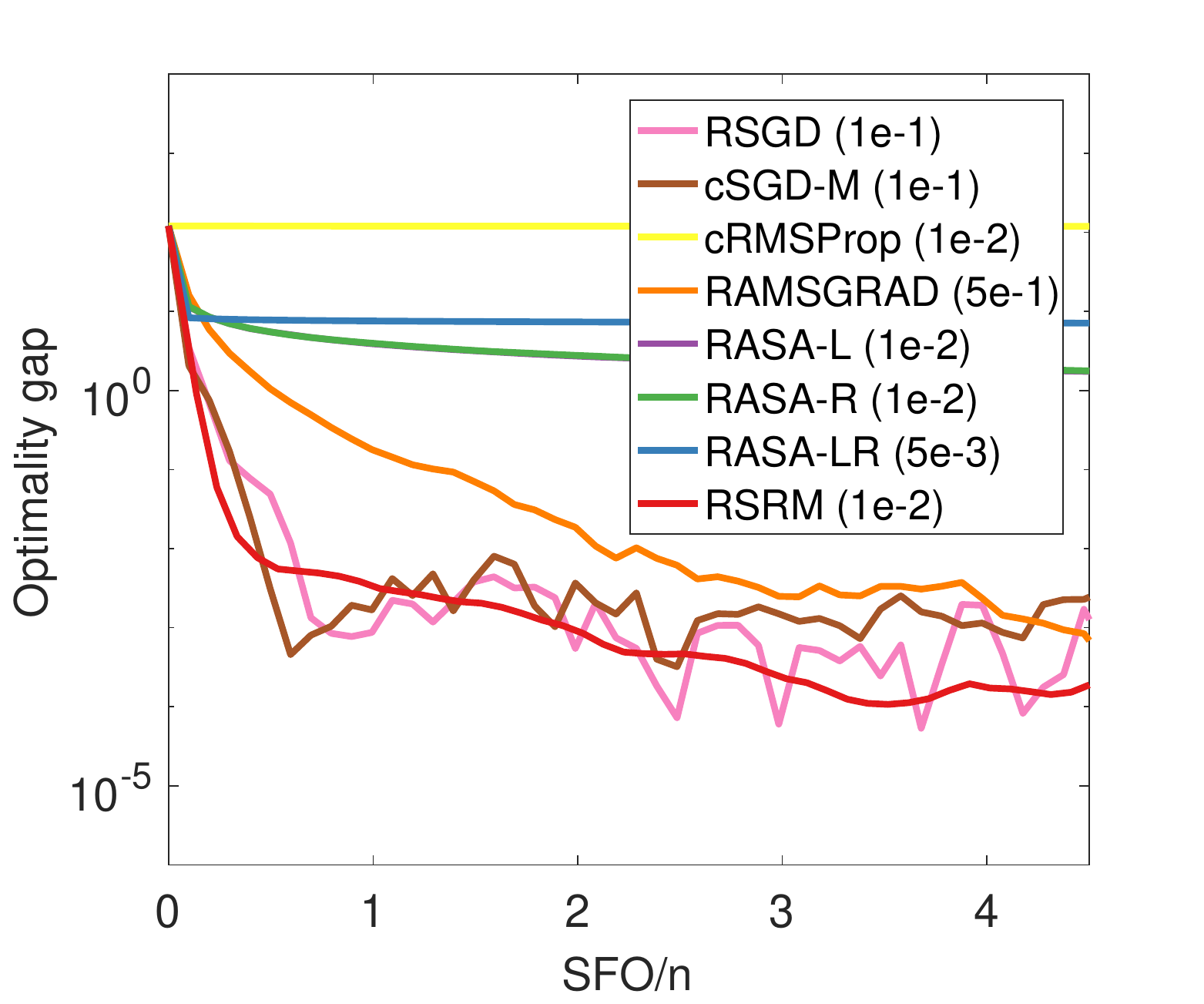}}
    \hspace{0.01in}
    \subfloat[Optimality gap vs. SFO (\texttt{KYLBERG})]{\includegraphics[width = 0.32\textwidth, height = 0.25\textwidth]{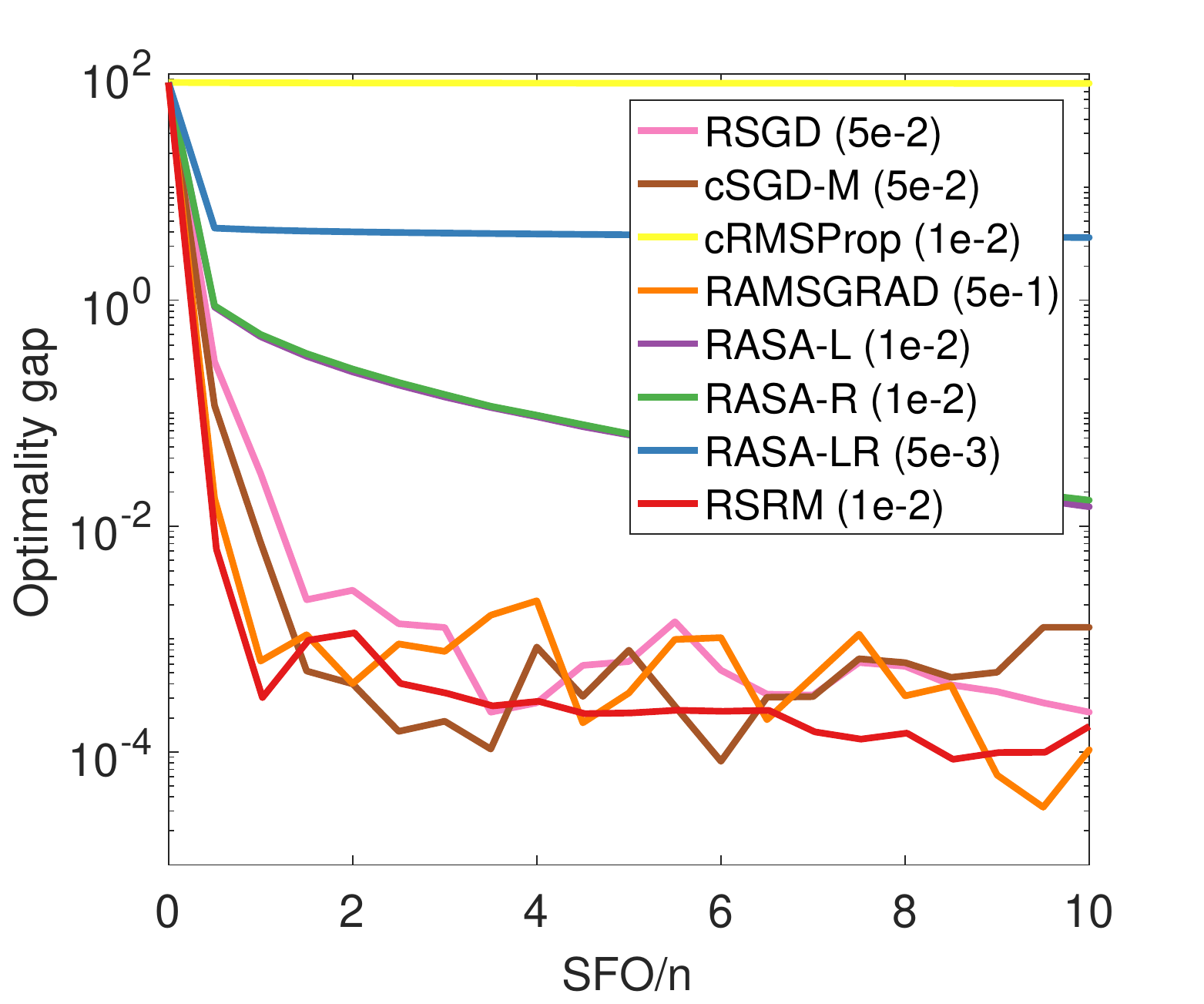}}\\
    \caption{RC problems on SPD manifold}
    \label{rc_spd}
\end{figure*}

All methods start with the same initialization and are terminated when the maximum iteration number is reached. For competing methods, we consider a square-root decaying step size $\eta_t = \eta_0/\sqrt{t}$, which is suggested in \cite{KasaiAdaptiveManifold2019}. We set the parameters of RSRM according to the theory, which is $\eta_t = \eta_0 / t^{1/3}$ and $\rho_t = \rho_0 / t^{2/3}$. A default value of $\rho_0 = 0.1$ provides good empirical performance. For all methods, $\eta_0$ are selected from $\{ 1, 0.5, 0.1, ..., 0.005, 0.001 \}$. The gradient momentum parameter in cSGD-M and RAMSGRAD is set to be $0.999$ and the adaptation momentum parameter in cRMSProp, RAMSGRAD and RASA is set to be $0.9$. We choose a mini-batch size of $5$ for RSRM and $10$ for algorithms excluding RSRM. Hence, the per-iteration cost of gradient evaluation is identical across methods. The initial batch size for RSRM is fixed to be $100$ (except for the problem of ICA where it is set to be $200$). All algorithms are coded in Matlab compatible to ManOpt toolbox \cite{BoumalManOpt2014} and results are reported on a i5-8600 3.1GHz CPU processor.

We consider principal component analysis (PCA) on Grassmann manifold, joint diagonalization of independent component analysis (ICA) on Stiefel manifold and computing Riemannian centroid (RC) on SPD manifold. Stiefel manifold $\text{St}(r, d) = \{ \mathbf X \in \mathbb{R}^{d \times r} : \mathbf X^T \mathbf{X} = \mathbf I_r \}$ is defined as the set of $d \times r$ column orthonormal matrices, which is a natural embedded submanifold of $\mathbb{R}^{d \times r}$. Therefore orthogonal projection is a valid vector transport. Grassmann manifold $\mathcal{G}(r, d)$ is the set of $r$-dimensional subspaces in $\mathbb{R}^d$. One representation of Grassmann manifold is by a matrix $\mathbf X \in \mathbb{R}^{d \times r} $ with orthonormal columns that span the subspace. This representation is not unique. Indeed, any $\mathbf {XR}$ for $\mathbf R \in O(r)$ is equivalent to $\mathbf X$, where $O(r)$ is the orthogonal group of dimension $r$. Hence, Grassmann manifold can be viewed as a quotient of Stiefel manifold, written as $\text{St}(r,d)/O(r)$. Finally, SPD manifold $\mathcal{S}_{++}^d$ is the set of $d \times d$ symmetric positive definite matrices, which forms the interior of a convex cone embedded in $\mathbb{R}^{n(n+1)/2}$.

\subsection{PCA on Grassmann manifold}

Consider a set of $n$ samples, represented by $\mathbf x_i \in \mathbb{R}^d, i = 1,...,n$. PCA aims to find a subspace where projection onto this subspace minimizes reconstruction error. This naturally defines a problem on Grassmann manifold, written as $\min_{\mathbf U \in \mathcal{G}(r, d)} \frac{1}{n} \sum_{i=1}^n \| \mathbf x_i - \mathbf{UU}^T \mathbf x_i \|^2$. This can be further simplified into  $\min_{\mathbf U \in \mathcal{G}(r, d)} \frac{1}{n} \sum_{i=1}^n \mathbf x_i^T \mathbf{UU}^T \mathbf x_i $. We first test RSRM on a baseline synthetic dataset (\texttt{SYN1}) with $(n, d, r) = (10^4, 10^2, 10)$. Then we increase dimension to $d = 500$ (\texttt{SYN2}) and consider a higher rank case with $r = 20$ (\texttt{SYN3}). In addition, we also consider two empirical datasets, \textsc{mnist} \cite{LecunMNIST1998} with $(n, d, r) = (60000, 784, 10)$ and \textsc{covtype} from LibSVM \cite{ChangLIBSVM2011} with $(n, d, r) = (581012, 54, 10)$. We measure the performance in optimality gap, which calculates the difference between current function value to the minimum, pre-calculated using Matlab function \textsc{pca}. Convergence results and the best-tuned $\eta_0$ are shown in Figure \ref{pca_grassmann}. We find that RSRM consistently outperforms others on both synthetic datasets (Figure \ref{pca_grassmann}\subref{sfo_syn1} to \subref{sfo_syn3}) as well as on real datasets (Figure \ref{pca_grassmann}\subref{sfo_mnist} and \subref{sfo_covtype}). It is also observed that on `easy' datasets, such as \texttt{SYN1} and \texttt{SYN2}, adaptive methods perform similarly compared with well-tuned SGD. Figure \ref{pca_grassmann}\subref{time_syn1} further illustrates the convergence behaviour against runtime where RSRM still maintains its superiority due to high per-iteration efficiency. Indeed, compared to SGD, RSRM only needs one extra vector transport operation, which can be as efficient as orthogonal projection for both Grassmann and Stiefel manifold. Other methods, such as cRMSProp, may further require large matrix computations for adaptation. In addition, Figure \ref{pca_grassmann}\subref{eta_sensitivity} and \subref{b_sensitivity} illustrates how convergence performance of RSRM varies under different choices of initial step size and batch size. It is noticed that RSRM with $\eta_0 = 0.005, 0.001$ outperforms best-tuned SGD. Also, RSRM seems to be insensitive to the initial batch size and surprisingly, larger batch size provides no benefit for its convergence under this problem.

\subsection{ICA on Stiefel manifold}
Independent component analysis (or blind source separation) aims to recover underlying components of observed multivariate data by assuming mutual independence of source signals. Joint diagonalization is a useful pre-processing step that searches for a pseudo-orthogonal matrix (i.e. Stiefel matrix) \cite{TheisICAJD2009} by solving $\min_{\mathbf U \in \text{St}(r,d)} - \frac{1}{n} \sum_{i=1}^n \|  \text{diag}(\mathbf U^T \mathbf X_i \mathbf U) \|^2_F$ with diag$(\mathbf A)$ returning diagonal elements of matrix $\mathbf A$. The symmetric matrices $\mathbf X_i \in \mathbb{R}^{d \times d}$ can be time-lagged covariance matrices \cite{BelouchraniTLMATRIX1997} or cumulant matrices \cite{CardosoCUMATRIX1999} constructed from the observed signals. We consider three image datasets described as follows. \texttt{YALEB} \cite{WrightYALEB2008} collects $n = 2414$ face images taken from various lighting environments. \texttt{CIFAR100} \cite{KrizhevskyCIFAR1002009} contains $n = 60000$ images of $100$ objects and \texttt{COIL100} \cite{NeneCOIL1001996} is made up of $n = 7200$ images from $100$ classes. To construct covariance representations from these datasets, we first downsize each image to $32 \times 32$ before applying Gabor-based kernel to extract Gabor features. Then the feature information is fused in a region covariance descriptors of size $43 \times 43$. We choose $r = d = 43$ for all problems and the results are presented in Figure \ref{ica_stiefel}. The optimal solution is obtained by running RSRM for sufficiently long. Similarly, we find that RSRM, although showing slow progress at the initial epochs, quickly converges to a lower value compared with others. This is mainly attributed to its variance reduction nature.

\subsection{RC on SPD manifold}
Computing Riemannian centroid on SPD manifold $\mathcal{S}_{++}^d$ are fundamental in many computer vision tasks, including object detection \cite{JayasumanaKRBF2015}, texture classification \cite{FarakiTCSPD2015} and particularly medical imaging \cite{ChengDTIP2012}. The problem concerns finding a mean representation of a set of SPD matrices, $\mathbf X_i \in \mathbb{R}^{d \times d}$. Among the many Riemannian metrics, Affine Invariant Riemannian Metric (AIRM) is most widely used to measure the closeness between SPD matrices. This induces a geodesic distance on SPD manifold, given by $d^2(\mathbf X_1, \mathbf X_2) = \| \log(\mathbf X_1^{-1/2} X_2 \mathbf X_1^{-1/2}) \|^2_F$ where $\log(\cdot)$ is the principal matrix logarithm. Riemannian centroid with respect to this distance is obtained by solving $\min_{\mathbf C \in \mathcal{S}_{++}^d} \frac{1}{n}\sum_{i=1}^n d^2(\mathbf C, \mathbf X_i)$. We first consider a simulated dataset consisting of $n = 5000$ SPD matrices in $\mathbb{R}^{10 \times 10}$, each with a condition number of $20$. Then we test our methods on \texttt{YALEB} face dataset \cite{WrightYALEB2008} and \texttt{KYLBERG} \cite{Kylbergkylberg2014} dataset that consists of $n = 4480$ texture images of $28$ classes. For each pixel, we generate a $5$-dimensional feature vector ($d = 5$), including pixel intensity, first- and second-order gradients. Subsequently, the covariance representation is similarly constructed for each image. Convergence results are shown in Figure \ref{rc_spd} where the optimal solutions are calculated by Riemannian Barzilai-Borwein algorithm \cite{IannazzoRBB2018}. By examining the figures, we also verify superiority of RSRM where it enjoys a more stable convergence due to variance reduction and sometimes converges to a lower objective value, as shown in Figure \ref{rc_spd}\subref{sfo_syn_rc} and \subref{sfo_yaleb_rc}. For the two real datasets, cRMSProp and RASA fails to perform comparably.

\section{Conclusion}
In this paper, we develop a one-sample online method that achieves the state-of-the-art lower bound complexity up to a logarithmic factor. This improves on SGD-based adaptive methods by a factor of $\Tilde{\mathcal{O}}(\epsilon^{-1})$. In particular, we use the stochastic recursive momentum estimator that only requires $\mathcal{O}(1)$ per-iteration gradient computation and achieves variance reduction without restarting the algorithm with a large batch gradient. Our experiment findings confirm the superiority of our proposed algorithm.

\bibliographystyle{aaai}

%\newpage
\onecolumn

\par\noindent\rule{\textwidth}{1pt}
\begin{center}
    \Large \textbf{Supplementary Material}
\end{center}
\par\noindent\rule{\textwidth}{1pt}

\appendix
\section{Proof of Lemma 1}

\begin{proof}
By definition of $\mathcal{F}_t$, we have $\mathbb{E}\| d_{t} - \text{grad}F(x_t) \|^2 = \mathbb{E}[\mathbb{E}[ \| d_{t} - \text{grad}F(x_t) \|^2 | \mathcal{F}_t ]]$. Then
\begin{align}
    &\mathbb{E}[ \| d_t - \text{grad}F(x_t) \|^2 |\mathcal{F}_t ] \nonumber\\
    &= \mathbb{E}[ \| (1 - \rho_t) \mathcal{T}_{x_{t-1}}^{x_t} (d_{t-1} - \text{grad}f_{\mathcal{S}_t}(x_{t-1})) + \text{grad}f_{\mathcal{S}_t}(x_t) - \text{grad}F(x_t) \|^2 | \mathcal{F}_t ] \nonumber\\
    &= \mathbb{E} [ \| (1 - \rho_t) (\mathcal{T}_{x_{t-1}}^{x_t} d_{t-1} + \text{grad}f_{\mathcal{S}_t}(x_t) - \mathcal{T}_{x_{t-1}}^{x_t} \text{grad}f_{\mathcal{S}_t}(x_{t-1}) ) + \rho_t \text{grad}f_{\mathcal{S}_t}(x_t) - \text{grad}F(x_t) \|^2  | \mathcal{F}_t] \nonumber\\
    &= \mathbb{E}[ \| (1-\rho_t) \mathcal{T}_{x_{t-1}}^{x_t} (d_{t-1} - \text{grad}F(x_{t-1})) + \rho_t (\text{grad}f_{\mathcal{S}_t}(x_t)  - \text{grad}F(x_t) ) \nonumber\\
    &+ (1-\rho_t) (  \text{grad}f_{\mathcal{S}_t}(x_t) - \mathcal{T}_{x_{t-1}}^{x_t} \text{grad}f_{\mathcal{S}_t}(x_{t-1}) +  \mathcal{T}_{x_{t-1}}^{x_t} \text{grad}F(x_{t-1}) - \text{grad}F(x_t) ) \|^2 | \mathcal{F}_t] \nonumber\\
    &= (1-\rho_t)^2 \| d_{t-1} - \text{grad}F(x_{t-1}) \|^2 + \mathbb{E}[ \| \rho_t (\text{grad}f_{\mathcal{S}_t}(x_t)  - \text{grad}F(x_t) ) \nonumber\\
    &+ (1-\rho_t) (  \text{grad}f_{\mathcal{S}_t}(x_t) - \mathcal{T}_{x_{t-1}}^{x_t} \text{grad}f_{\mathcal{S}_t}(x_{t-1}) +  \mathcal{T}_{x_{t-1}}^{x_t} \text{grad}F(x_{t-1}) - \text{grad}F(x_t) )   \|^2 | \mathcal{F}_t ] \nonumber\\
    &\leq (1-\rho_t)^2 \| d_{t-1} - \text{grad}F(x_{t-1}) \|^2 + 2 \rho_t^2 \mathbb{E}[\| \text{grad}f_{\mathcal{S}_t}(x_t)  - \text{grad}F(x_t)  \|^2 | \mathcal{F}_t] \nonumber\\
    &+ 2 (1- \rho_t)^2 \mathbb{E}[\| \text{grad}f_{\mathcal{S}_t}(x_t) - \mathcal{T}_{x_{t-1}}^{x_t} \text{grad}f_{\mathcal{S}_t}(x_{t-1}) +  \mathcal{T}_{x_{t-1}}^{x_t} \text{grad}F(x_{t-1}) - \text{grad}F(x_t) \|^2 |\mathcal{F}_t] \nonumber\\
    &\leq (1-\rho_t)^2 \| d_{t-1} - \text{grad}F(x_{t-1}) \|^2 + 2 \rho_t^2 \mathbb{E}[\| \text{grad}f_{\mathcal{S}_t}(x_t)  - \text{grad}F(x_t)  \|^2 | \mathcal{F}_t] \nonumber\\
    &+ 2(1 - \rho_t)^2 \mathbb{E}[\| \text{grad}f_{\mathcal{S}_t}(x_t) - \mathcal{T}_{x_{t-1}}^{x_t} \text{grad}f_{\mathcal{S}_t}(x_{t-1}) \|^2 |\mathcal{F}_t] \nonumber\\
    &= (1-\rho_t)^2 \| d_{t-1} - \text{grad}F(x_{t-1}) \|^2 + \frac{2\rho_t^2}{|\mathcal{S}_t|} \mathbb{E}_\omega \| \text{grad}f(x_t, \omega) - \text{grad}F(x_t)\|^2 + \frac{2(1-\rho_t)^2}{|\mathcal{S}_t|} \mathbb{E}_\omega \| \text{grad}f(x_t, \omega) \nonumber\\
    &- \mathcal{T}_{x_{t-1}}^{x_t} \text{grad}f(x_{t-1}, \omega)  \|^2 \nonumber\\
    &\leq (1-\rho_t)^2 \| d_{t-1} - \text{grad}F(x_{t-1}) \|^2 + \frac{2\rho_t^2 \sigma^2}{|\mathcal{S}_t|} + \frac{2(1-\rho_t)^2 \eta_{t-1}^2 \Tilde{L}^2}{|\mathcal{S}_t|} \| d_{t-1}\|^2 \nonumber\\
    &\leq (1-\rho_t)^2 \big( 1 + \frac{4\eta_{t-1}^2 \Tilde{L}^2}{|\mathcal{S}_t|} \big) \| d_{t-1} - \text{grad}F(x_{t-1}) \|^2 + \frac{4(1-\rho_t)^2 \eta_{t-1}^2 \Tilde{L}^2}{|\mathcal{S}_t|} \| \text{grad}F(x_{t-1}) \|^2 + \frac{2\rho_t^2 \sigma^2}{|\mathcal{S}_t|}. 
\end{align}
where the last equality uses isometry property of vector transport $\mathcal{T}_{x_{t-1}}^{x_t}$ and the fact that it is measurable in $\mathcal{F}_t$. Also, we use the unbiasedness of stochastic gradient $\text{grad}f_{\mathcal{S}_t}(x_t)$. The first and last inequalities follow from $\| a + b \|^2 \leq 2\| a\|^2 + 2 \| b\|^2$. The second inequality applies $\mathbb{E}\|x - \mathbb{E}[x] \|^2 \leq \mathbb{E}\|x\|^2$ for random variable $x$. The second last inequality follows from Assumptions \ref{BiasVarAssump} and \ref{MSLipAssump}. By taking full expectation, we obtain the desired result.

\end{proof}

\section{Proof of Theorem 1}

\begin{proof}
By retraction $L$ smoothness of $F$, we have
\begin{align}
    F(x_{t+1}) &\leq F(x_t) - \langle \text{grad}F(x_t), \eta_t d_t \rangle + \frac{\eta_t^2 L}{2} \| d_t \|^2  \nonumber\\
    &= F(x_t) - \frac{\eta_t}{2} \| \text{grad}F(x_t) \|^2 - \frac{\eta_t}{2} \|d_t \|^2 + \frac{\eta_t}{2} \| d_t - \text{grad}F(x_t) \|^2 + \frac{L \eta_t^2}{2} \| d_t\|^2 \nonumber\\
    &= F(x_t) - \frac{\eta_t}{2} \| \text{grad}F(x_t) \|^2 + \frac{\eta_t}{2} \| d_t - \text{grad}F(x_t) \|^2 - (\frac{\eta_t}{2} - \frac{L\eta_t^2}{2}) \|d_t \|^2 \nonumber\\
    &\leq F(x_t) - \frac{\eta_t}{2} \| \text{grad}F(x_t) \|^2 + \frac{\eta_t}{2} \| d_t - \text{grad}F(x_t) \|^2,
\end{align}
where for the last inequality we choose $\eta_t \leq \frac{1}{L}$.Given $\eta_t = c_\eta (t+1)^{-1/3}$, we can ensure this condition by requiring $c_\eta \leq \frac{1}{L}$. Given the choice that $\eta_t = c_\eta (t+1)^{-1/3}$, $\rho_t = c_\rho t^{-2/3}$, $|\mathcal{S}_t| = b$ for all $t$, we construct a Lyapunov function 
\begin{equation}
    R_t := \mathbb{E}[F(x_t)] + \frac{C}{\eta_{t-1}} {\mathbb{E}\| d_t - \text{grad}F(x_t) \|^2},
\end{equation}
for some constant $C > 0$. Denote estimation error as $s_t = d_t - \text{grad}F(x_t)$. Then we can bound the scaled difference between two consecutive estimation error as
\begin{align}
    \frac{\mathbb{E}\| s_{t+1} \|^2}{\eta_{t}} - \frac{\mathbb{E}\| s_t \|^2}{\eta_{t-1}} &\leq \frac{(1-\rho_{t+1})^2 \big( 1 + \frac{4\eta_{t}^2 \Tilde{L}^2}{b} \big) \mathbb{E}\| s_t \|^2 + \frac{4(1-\rho_{t+1})^2 \eta_{t}^2 \Tilde{L}^2}{b} \mathbb{E}\| \text{grad}F(x_{t}) \|^2 + \frac{2\rho_{t+1}^2 \sigma^2}{b}}{\eta_{t}} - \frac{\mathbb{E}\| s_t \|^2}{\eta_{t-1}} \nonumber\\
    &= \big(\frac{1}{\eta_{t}} (1 - \rho_{t+1})^2 (1+\frac{4\eta_t^2 \Tilde{L}^2}{b}) - \frac{1}{\eta_{t-1}} \big) \mathbb{E}\| s_t \|^2 + \frac{4(1 - \rho_{t+1})^2 \eta_t \Tilde{L}^2}{b } \mathbb{E}\| \text{grad}F(x_t) \|^2 + \frac{2\rho_{t+1}^2 \sigma^2}{b \eta_{t}} \nonumber\\
    &\leq \big( \frac{1}{\eta_{t}} - \frac{1}{\eta_{t-1}} + \frac{4\eta_t \Tilde{L}^2}{b} - \frac{\rho_{t+1}}{\eta_t} \big) \mathbb{E}\| s_t \|^2 + \frac{4 \eta_t \Tilde{L}^2}{b } \mathbb{E}\| \text{grad}F(x_t) \|^2 + \frac{2 \sigma^2}{b (t+1)}, \label{ruiurtiti}
\end{align}
where the second inequality uses the fact that $(1-\rho_{t+1})^2 \leq 1- \rho_{t+1} \leq 1$. Next we bound the first and third terms as follows. The third term $\sum_{t=1}^{T} \frac{2\sigma^2}{b(t+1)} \leq \frac{2\sigma^2 \ln(T+1)}{b}$ from the bound on Harmonic series. Now we bound the first term. Consider the convex function $h(x) := x^{1/3}$. By first order characterization, $h(x+1) \leq h(x) + h'(x) = x^{1/3} + \frac{1}{3}x^{-2/3}$. Therefore, we have $1/\eta_t - 1/\eta_{t-1} \leq \frac{c_\eta}{3} t^{-2/3} \leq \frac{c_\eta}{3} (t+1)^{-1/3} = \frac{\eta_t}{3}$, where we can easily verify that $t^{-2/3} \leq (t+1)^{-1/3}$ for $t \geq 2$.  Also given that $\frac{\rho_{t+1}}{\eta_t} = (\frac{10 \Tilde{L}^2}{b} + \frac{1}{3}) c_\eta (t+1)^{-1/3} = (\frac{10\Tilde{L}^2}{b} + \frac{1}{3}) \eta_t$ by the choice that $c_\rho/c_\eta = (\frac{ 10 \Tilde{L}^2}{b} + \frac{1}{3}) c_\eta$. Therefore, 
\begin{equation}
    \frac{1}{\eta_t} - \frac{1}{\eta_{t-1}} + \frac{4\eta_t \Tilde{L}^2}{b} - \frac{\rho_{t+1}}{\eta_t} \leq \frac{\eta_t}{3} + \frac{4\eta_t \Tilde{L}^2}{b} - (\frac{10\Tilde{L}^2}{b} + \frac{1}{3}) \eta_t = - \frac{6\Tilde{L}^2 \eta_t}{b}. 
\end{equation}
Then substituting this results in \eqref{ruiurtiti} gives 
\begin{equation}
    \sum_{t=1}^T  \big( \frac{\mathbb{E}\| s_{t+1} \|^2}{\eta_{t}} - \frac{\mathbb{E}\| s_t \|^2}{\eta_{t-1}} \big) \leq - \frac{6\Tilde{L}^2}{b} \sum_{t=1}^T \eta_t \mathbb{E}\| s_t \|^2 + \frac{4\Tilde{L}^2}{b} \sum_{t=1}^T \eta_t \mathbb{E} \| \text{grad}F(x_t) \|^2 + \frac{2\sigma^2 \ln(T+1)}{b}. 
\end{equation}
Now choose $C = \frac{b}{12\Tilde{L}^2}$. Then we have
\begin{align}
    R_{t+1} - R_t &= \mathbb{E}[F(x_{t+1}) - F(x_t)] + \mathbb{E}[ \frac{b}{12 \Tilde{L}^2 \eta_t} \mathbb{E}\| s_{t+1} \|^2 - \frac{b}{12 \Tilde{L}^2 \eta_{t-1}} \mathbb{E}\|s_t \|^2 ] \nonumber\\
    &\leq -\frac{\eta_t}{2} \mathbb{E}\| \text{grad}F(x_t) \|^2 + \frac{\eta_t}{2}\mathbb{E}\| s_t \|^2 + \mathbb{E}[ \frac{b}{12 \Tilde{L}^2 \eta_t} \mathbb{E}\| s_{t+1} \|^2 - \frac{b}{12 \Tilde{L}^2 \eta_{t-1}} \mathbb{E}\|s_t \|^2 ].
\end{align}
Telescoping this result from $t = 1,...,T$ gives 
\begin{align}
    R_{T+1} - R_1 &\leq -\sum_{t=1}^T \frac{\eta_t}{2} \mathbb{E}\| \text{grad}F(x_t) \|^2 + \sum_{t=1}^T \frac{\eta_t}{2} \mathbb{E}\| s_t \|^2 - \sum_{t=1}^T \frac{\eta_t}{2} \mathbb{E}\| s_t \|^2 + \sum_{t=1}^T \frac{\eta_t}{3} \mathbb{E}\| \text{grad}F(x_t) \|^2 + \frac{\sigma^2 \ln(T+1) }{6 \Tilde{L}^2} \nonumber\\
    &= - \sum_{t=1}^T \frac{\eta_t}{6} \mathbb{E}\| \text{grad}F(x_t) \|^2 + \frac{\sigma^2 \ln(T+1) }{6 \Tilde{L}^2}.  
\end{align}
Given $\eta_t$ is a decreasing sequence, we have $\sum_{t=1}^T \frac{\eta_t}{6} \mathbb{E}\| \text{grad}F(x_t) \|^2 \geq \frac{\eta_T}{6} \sum_{t=1}^T \mathbb{E}\| \text{grad}F(x_t) \|^2$. Therefore, we obtain 
\begin{align}
    \frac{1}{T} \sum_{t=1}^T \mathbb{E}\| \text{grad}F(x_t) \|^2 \leq \frac{6 (R_1 - R_{T+1}) + \frac{\sigma^2 \ln(T+1)}{\Tilde{L}^2} }{\eta_T T} &\leq \frac{6 F(x_1) - 6\mathbb{E}[F(x_{T+1})] + \frac{b}{2\Tilde{L}^2}\mathbb{E} \| d_1 - \text{grad}F(x_t) \|^2 + \frac{\sigma^2 \ln(T+1)}{\Tilde{L}^2} }{c_\eta (T+1)^{-1/3} T} \nonumber\\
    &\leq \frac{6 \Delta + \frac{\sigma^2}{2 \Tilde{L}^2} + \frac{\sigma^2 \ln(T+1)}{\Tilde{L}^2}}{c_\eta T} (T+1)^{1/3} \nonumber\\
    &\leq \frac{6 \Delta + \frac{\sigma^2}{2 \Tilde{L}^2} + \frac{\sigma^2 \ln(T+1)}{\Tilde{L}^2}}{c_\eta T} (T^{1/3} + 1) 
    = \frac{M}{T^{2/3}} + \frac{M}{T} = \Tilde{\mathcal{O}} ( \frac{1}{T^{2/3}}),
\end{align}
where we use the fact that $(a+b)^{1/3} \leq a^{1/3} + b^{1/3}$ and $M := (6 \Delta + \frac{\sigma^2}{2\Tilde{L}^2} + \frac{\sigma^2 \ln(T+1)}{\Tilde{L}^2})/c_\eta$. Hence to achieve $\epsilon$-approximate solution, we require $\mathbb{E}\| \text{grad}F(\Tilde{x}) \|^2 = \frac{1}{T} \sum_{t=1}^T \mathbb{E}\| \text{grad}F(x_t) \|^2 \leq \epsilon^2$. Hence requiring $\frac{M}{T^{2/3}} \leq \epsilon^2$ is sufficient for this purpose, which gives $\Tilde{\mathcal{O}}(\epsilon^{-3})$

\end{proof}

\section{Additional experiment results}

\subsection{Optimality gap against runtime}
\begin{figure*}
\captionsetup{justification=centering}
    \centering
    \subfloat[PCA (\texttt{SYN2})]{\includegraphics[width = 0.32\textwidth, height = 0.25\textwidth]{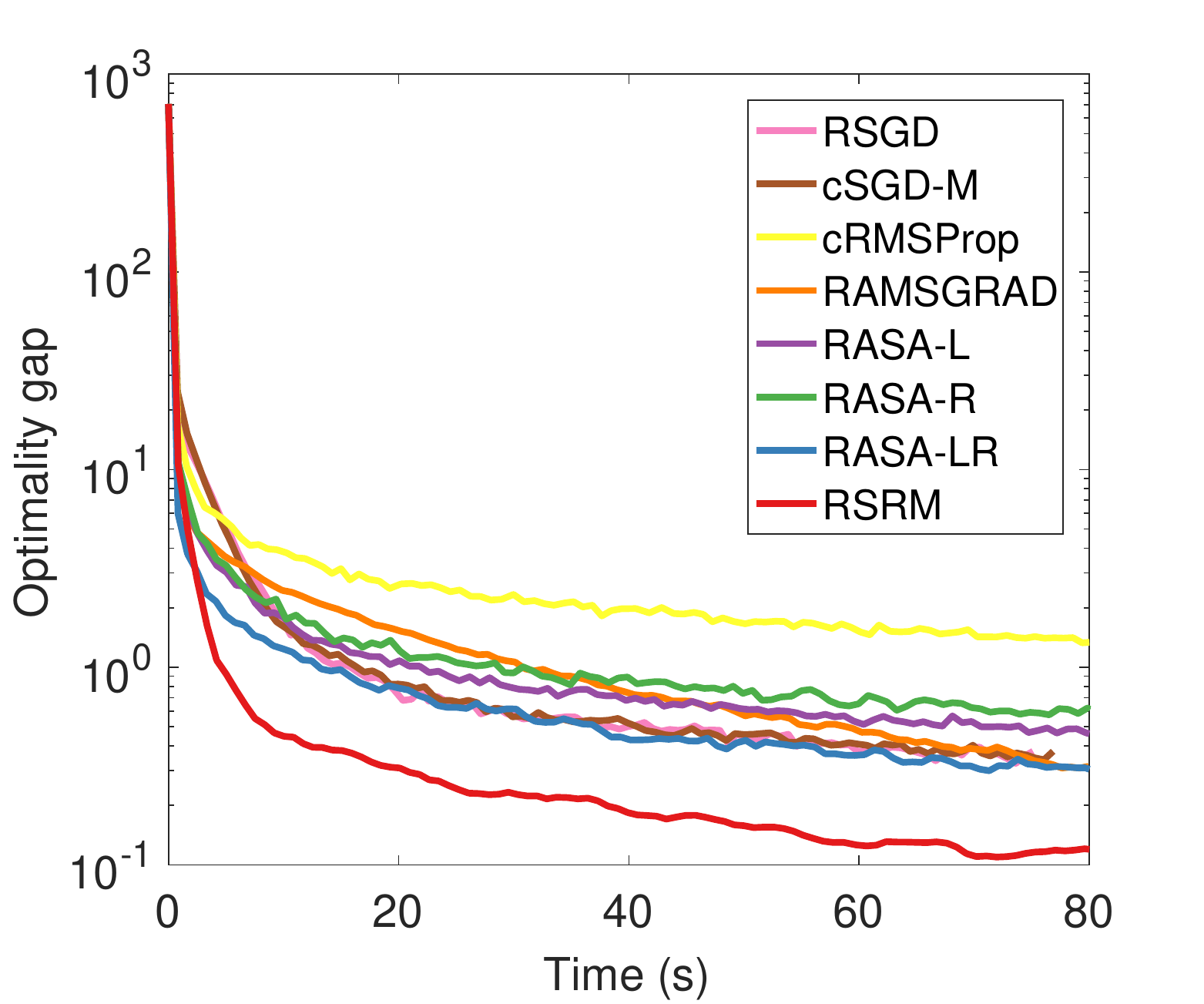}} 
    \hspace{0.01in}
    \subfloat[PCA (\texttt{SYN3})]{\includegraphics[width = 0.32\textwidth, height = 0.25\textwidth]{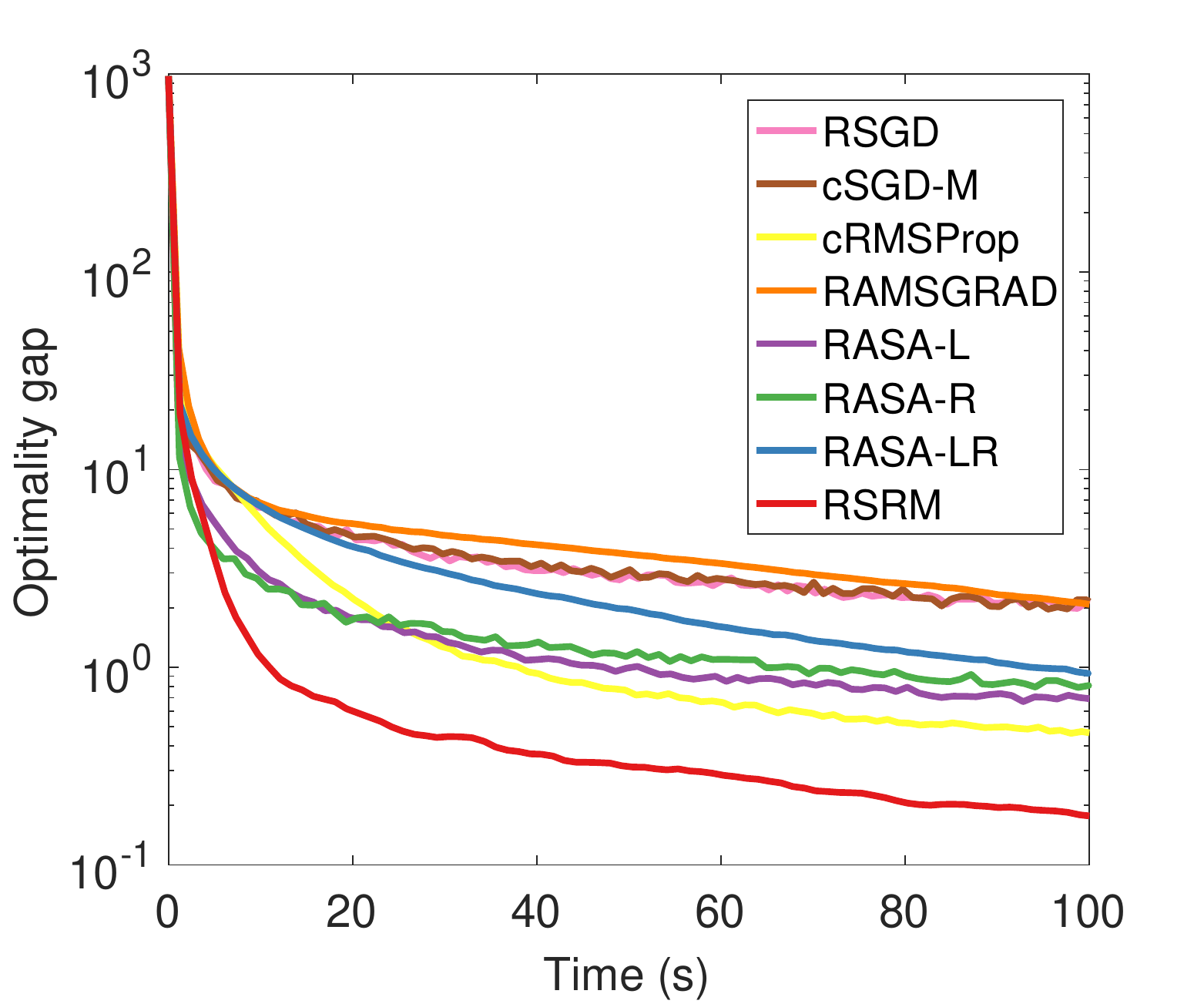}}
    \hspace{0.01in}
    \subfloat[PCA (\texttt{MNIST})]{\includegraphics[width = 0.32\textwidth, height = 0.25\textwidth]{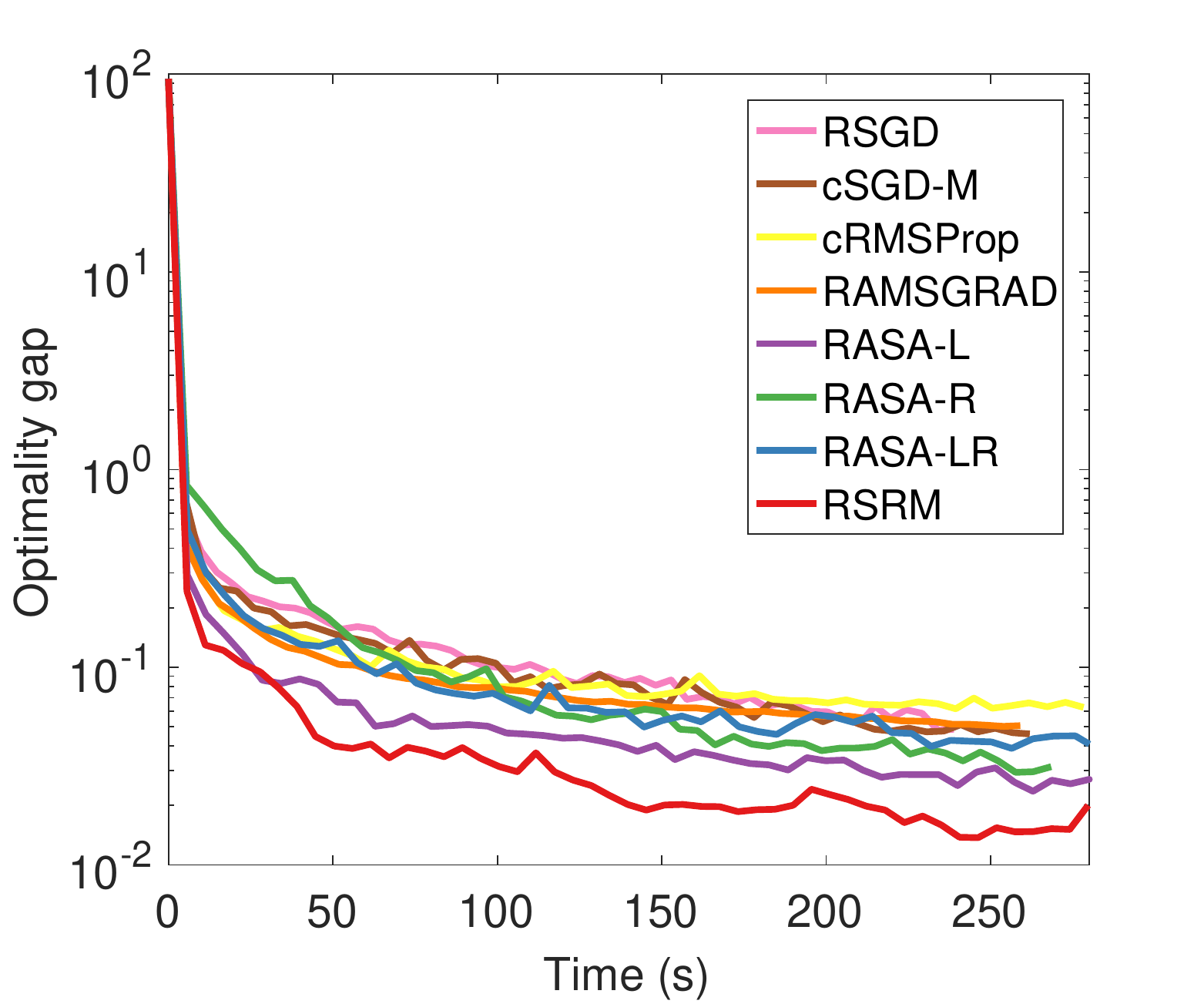}}\\[-0.2ex]
    \subfloat[PCA (\texttt{COVTYPE})]{\includegraphics[width = 0.32\textwidth, height = 0.25\textwidth]{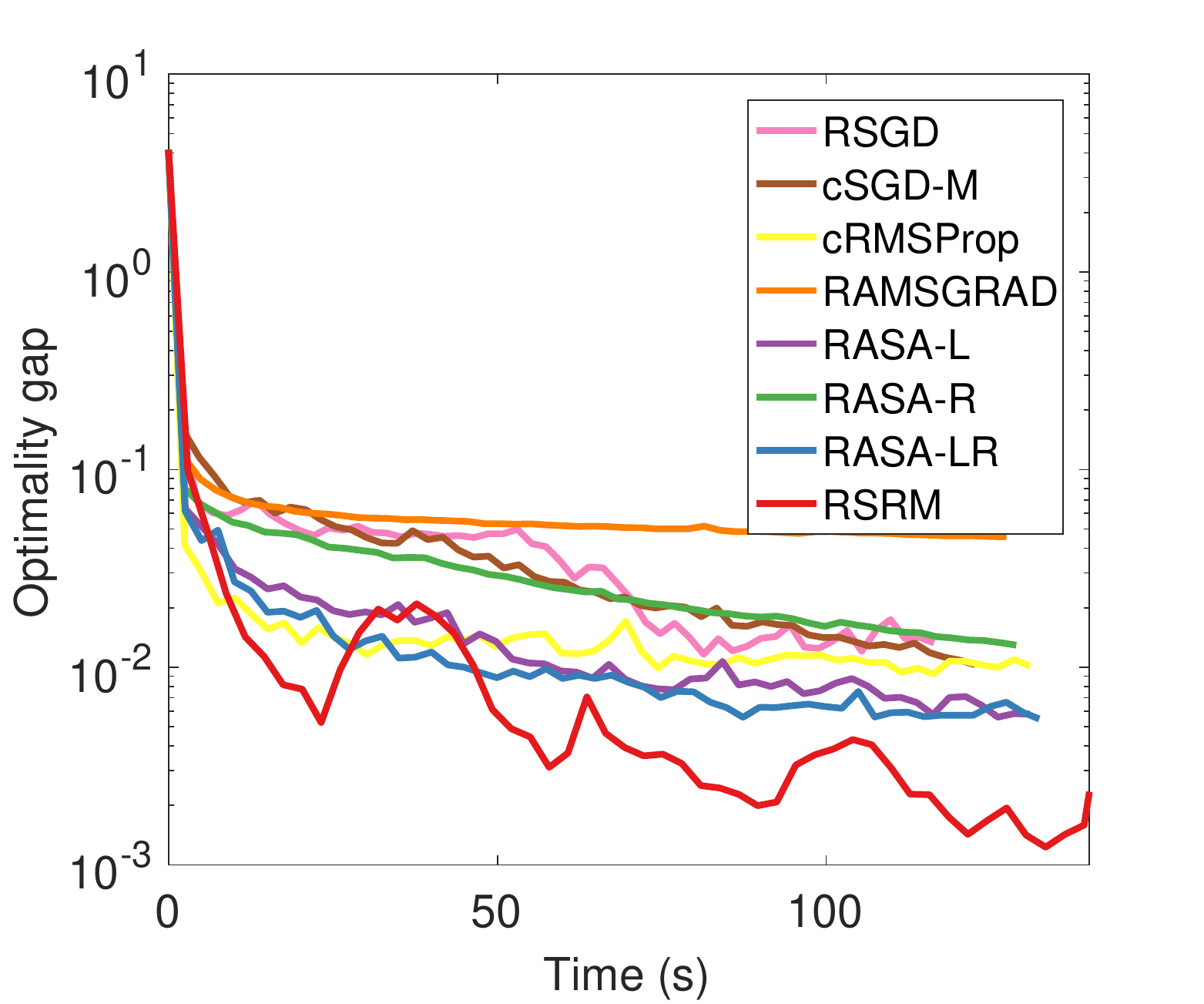}}
    \hspace{0.01in}
    \subfloat[ICA (\texttt{YALEB})]{\includegraphics[width = 0.32\textwidth, height = 0.25\textwidth]{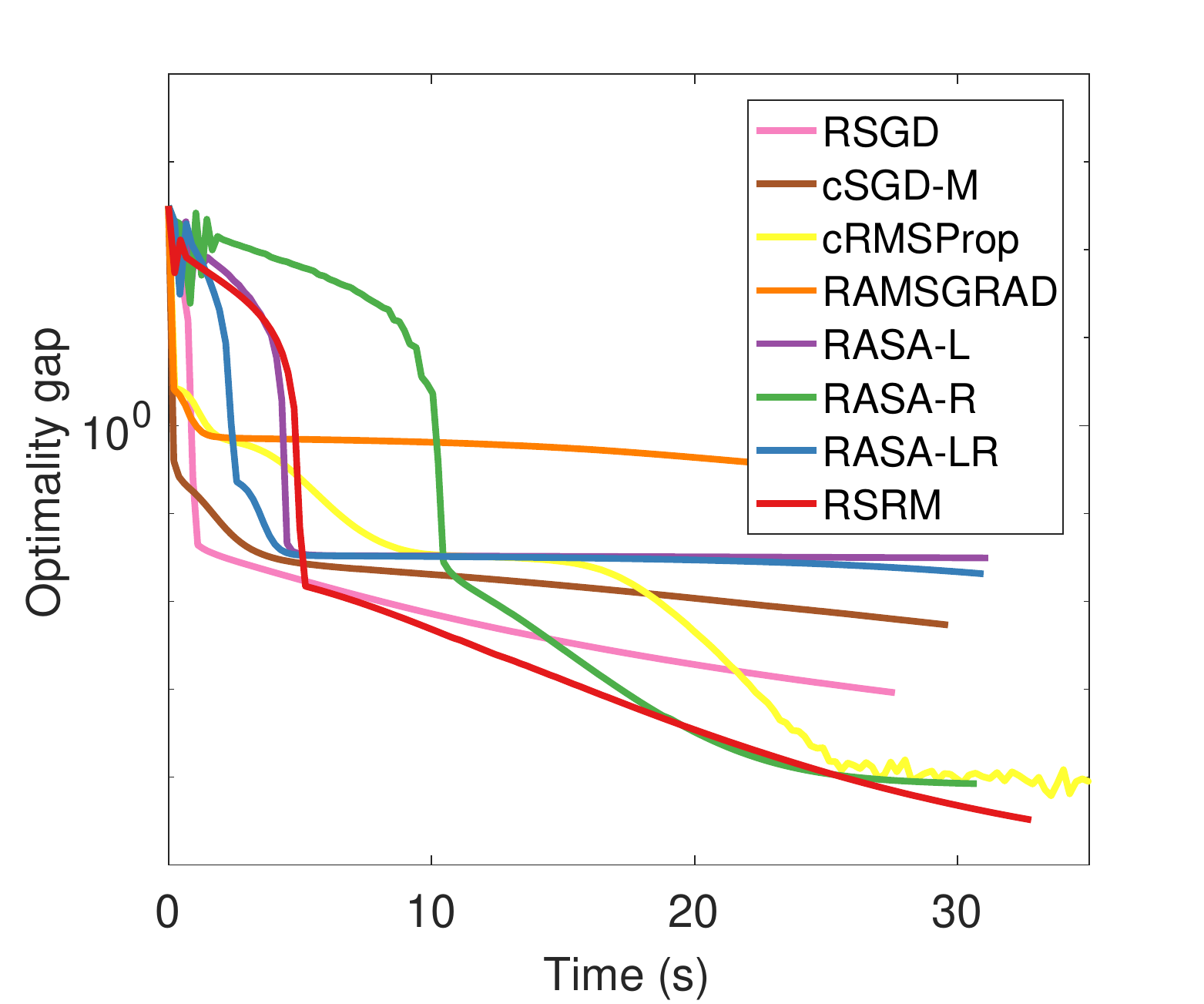}}
    \hspace{0.01in}
    \subfloat[ICA (\texttt{CIFAR100})]{\includegraphics[width = 0.32\textwidth, height = 0.25\textwidth]{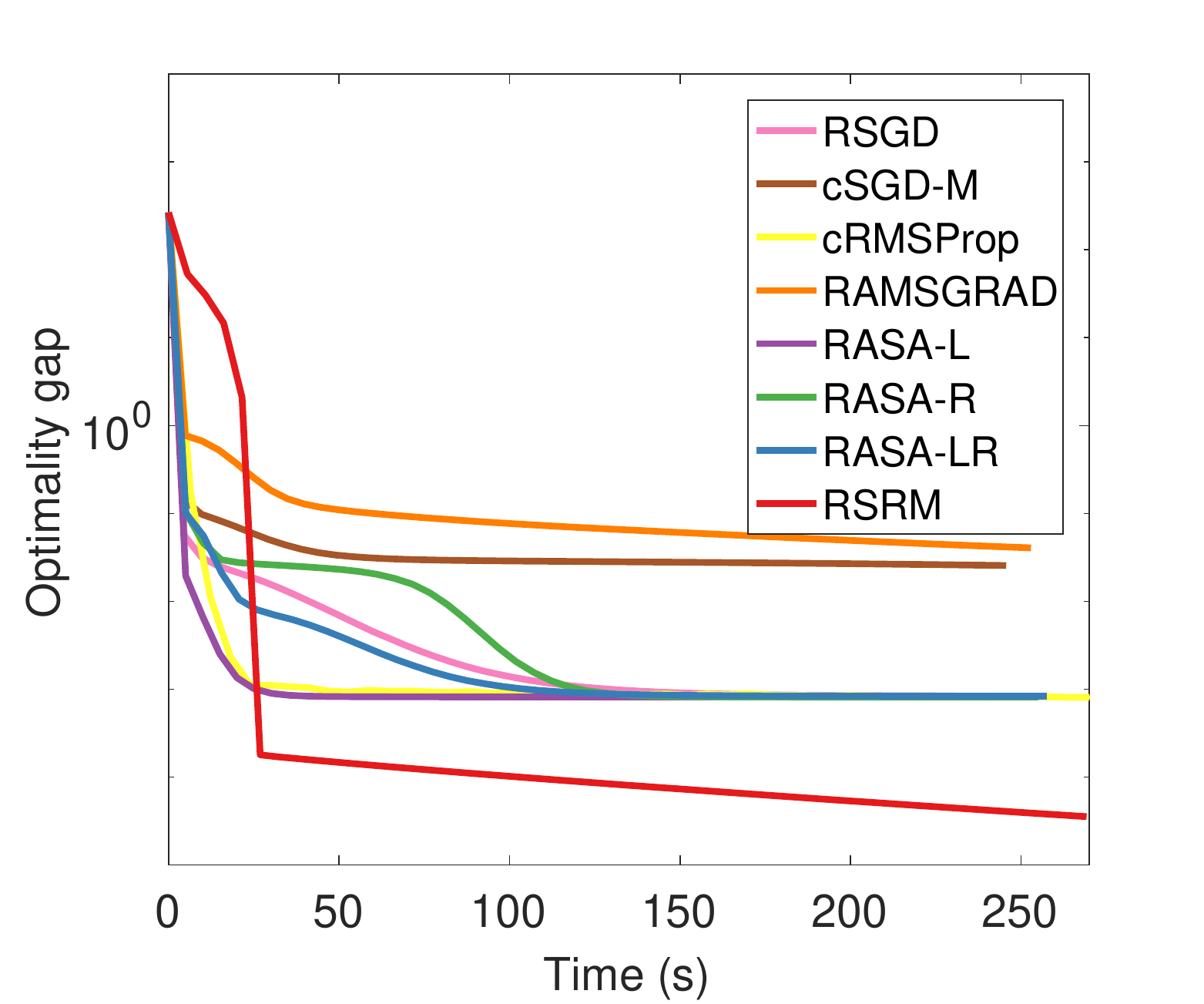}}\\[-0.2ex]
    \subfloat[ICA (\texttt{COIL100})]{\includegraphics[width = 0.32\textwidth, height = 0.25\textwidth]{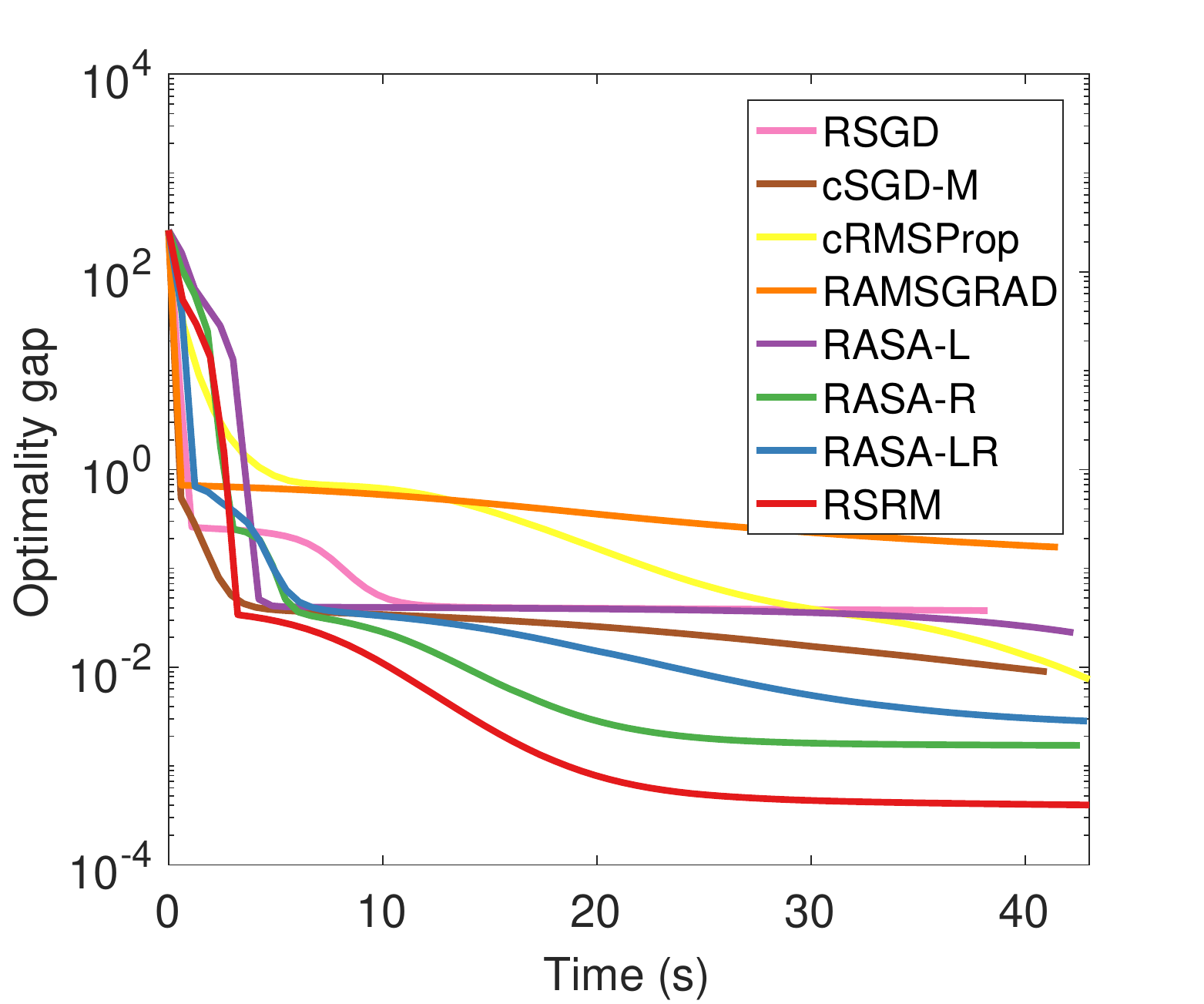}}
    \hspace{0.01in}
    \subfloat[RC (\texttt{SYN1})]{\includegraphics[width = 0.32\textwidth, height = 0.25\textwidth]{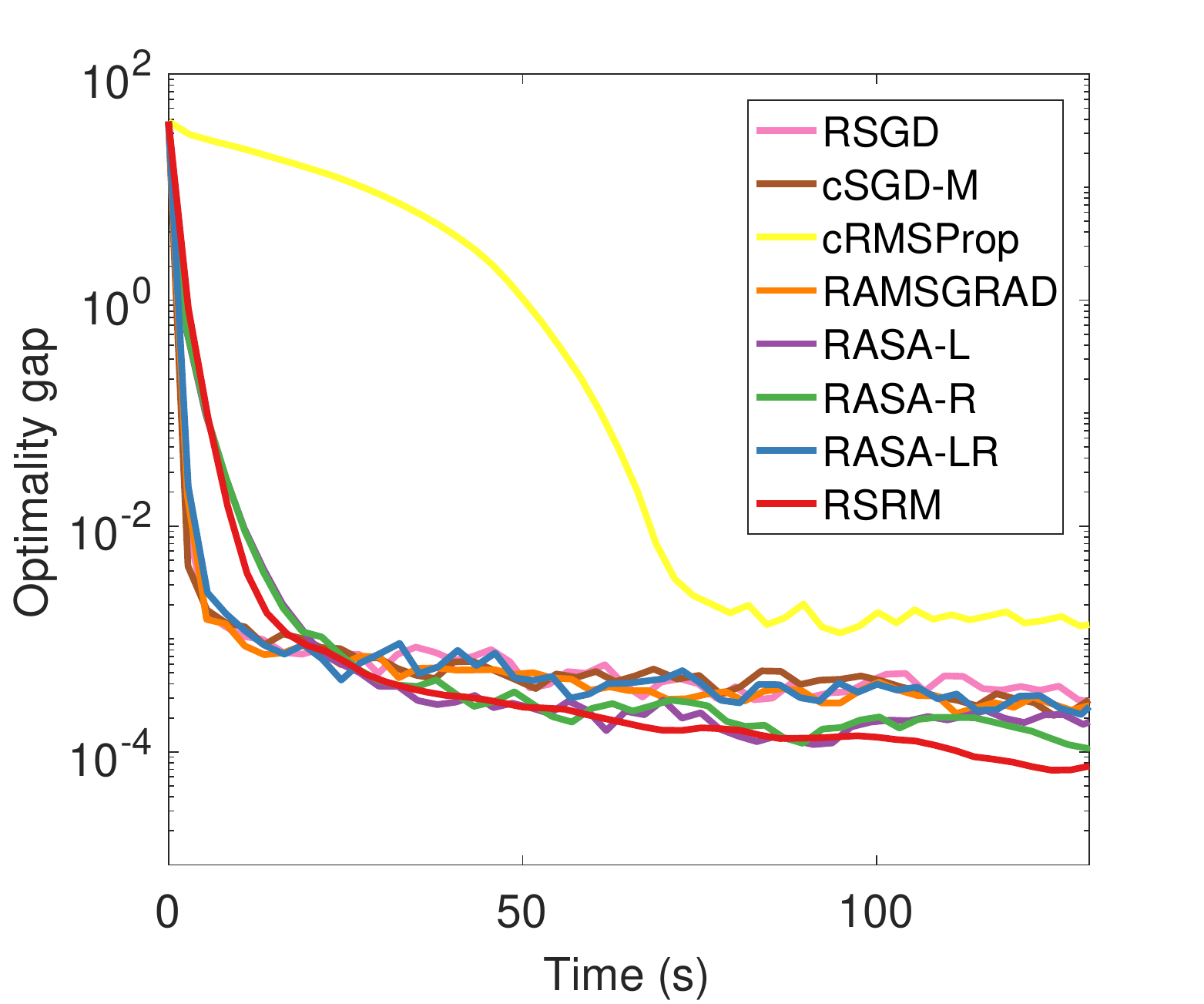}}
    \hspace{0.01in}
    \subfloat[RC (\texttt{YALEB})]{\includegraphics[width = 0.32\textwidth, height = 0.25\textwidth]{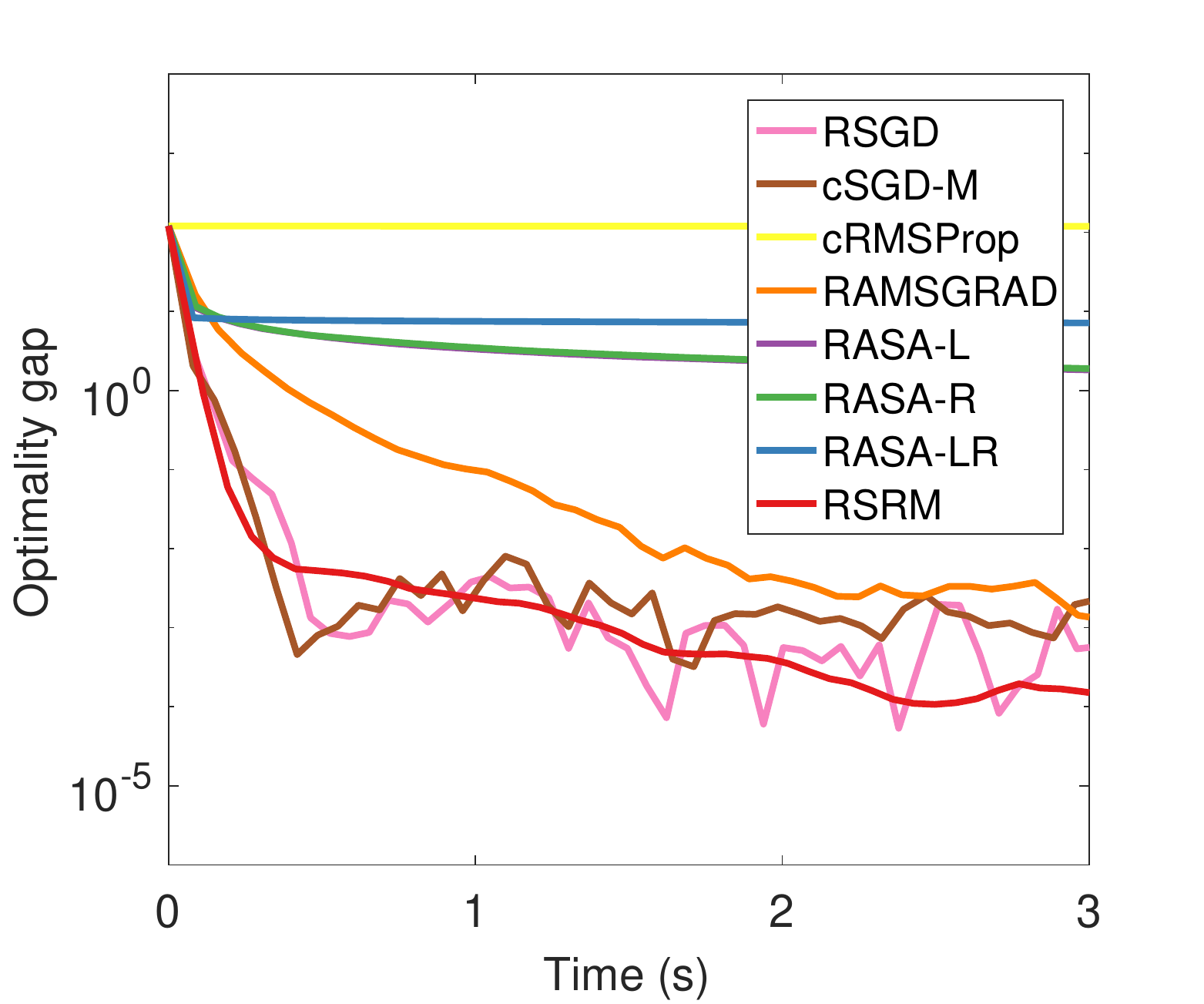}}\\[-0.2ex]
    \subfloat[RC (\texttt{KYLBERG})]{\includegraphics[width = 0.32\textwidth, height = 0.25\textwidth]{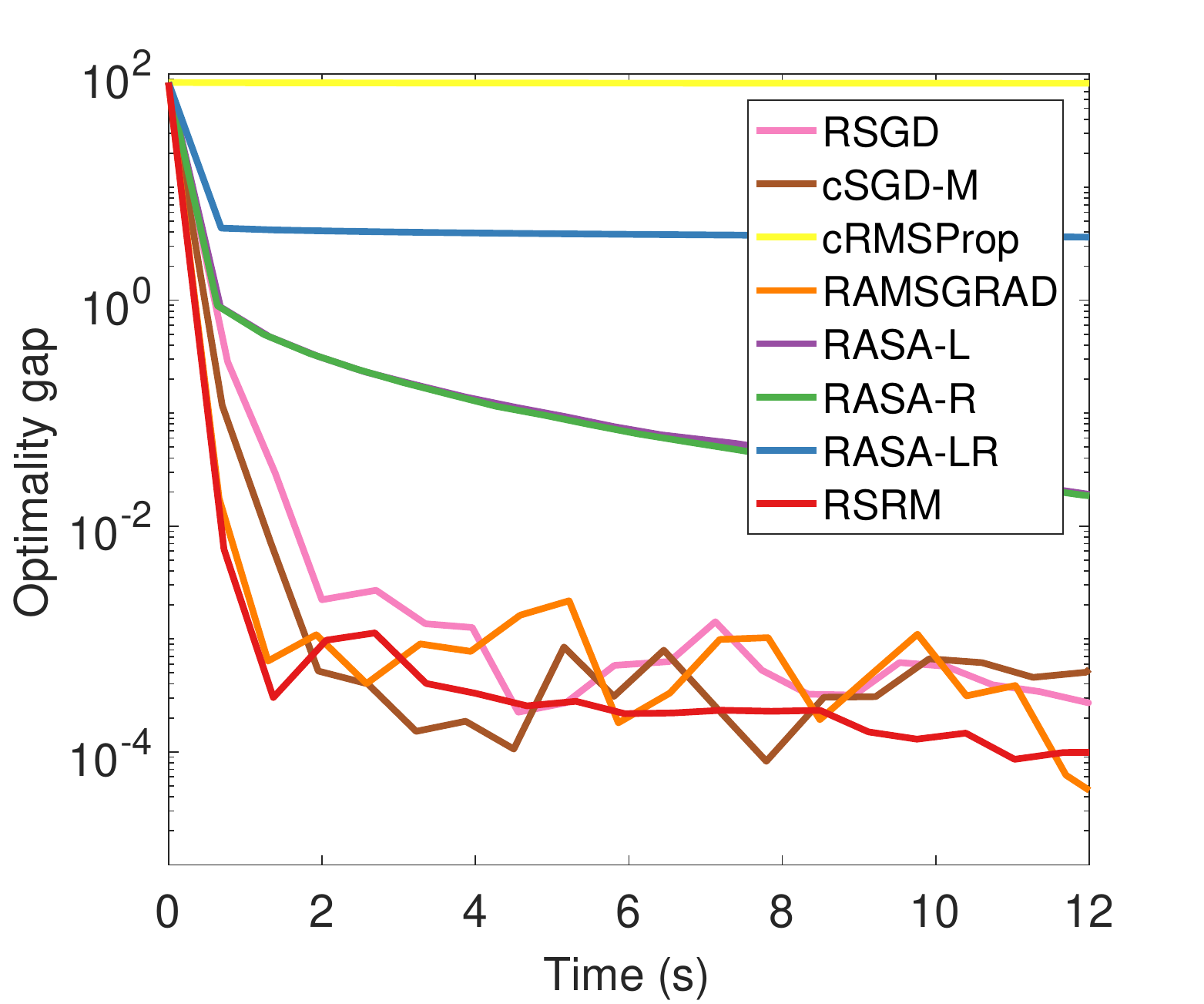}}
    \caption{Optimality gap against runtime for all problems} 
    \label{Runtime_all}
\end{figure*}

\end{document}